\documentclass[11pt,twoside]{article}
\usepackage{amsmath, amsthm, amscd, amsfonts, amssymb, graphicx, color}

\date{}

\begin{document}

\centerline{}

\centerline {\Large{\bf Atomic systems in $n$-Hilbert spaces and their tensor products }}

\newcommand{\mvec}[1]{\mbox{\bfseries\itshape #1}}
\centerline{}
\centerline{\textbf{Prasenjit Ghosh}}
\centerline{Department of Pure Mathematics, University of Calcutta,}
\centerline{35, Ballygunge Circular Road, Kolkata, 700019, West Bengal, India}
\centerline{e-mail: prasenjitpuremath@gmail.com}
\centerline{}
\centerline{\textbf{T. K. Samanta}}
\centerline{Department of Mathematics, Uluberia College,}
\centerline{Uluberia, Howrah, 711315,  West Bengal, India}
\centerline{e-mail: mumpu$_{-}$tapas5@yahoo.co.in}

\newtheorem{Theorem}{\quad Theorem}[section]

\newtheorem{definition}[Theorem]{\quad Definition}

\newtheorem{theorem}[Theorem]{\quad Theorem}

\newtheorem{remark}[Theorem]{\quad Remark}

\newtheorem{corollary}[Theorem]{\quad Corollary}

\newtheorem{note}[Theorem]{\quad Note}

\newtheorem{lemma}[Theorem]{\quad Lemma}

\newtheorem{example}[Theorem]{\quad Example}

\newtheorem{result}[Theorem]{\quad Result}
\newtheorem{conclusion}[Theorem]{\quad Conclusion}

\newtheorem{proposition}[Theorem]{\quad Proposition}

\begin{abstract}
\textbf{\emph{Concept of a family of local atoms in \,$n$-Hilbert space is being studied.\,$K$-frame in tensor product of \,$n$-Hilbert spaces is described and a characterization is given.\,Atomic system in tensor product of \,$n$-Hilbert spaces is presented and established a relationship between atomic systems in \,$n$-Hilbert spaces and their tensor products.}}
\end{abstract}
{\bf Keywords:}  \emph{Atomic system, $K$-frame, Tensor product of Hilbert spaces, linear\\ \smallskip\hspace{2.5cm}$n$-normed space, $n$-Hilbert space.}

\section{Introduction}

\smallskip\hspace{.6 cm}\,In recent times, many generalizations of frames have been appeared.\,Some of them are \,$K$-frame, \,$g$-frame, fusion frame and \,$g$-fusion frame etc.\;$K$-frames for a separable Hilbert spaces were introduced by Lara Gavruta \cite{L} to study atomic decomposition systems for a bounded linear operator.\,Infact, generalized atomic subspaces for operators in Hilbert spaces were studied by P.\,Ghosh and T.\,K.\,Samanta {\cite{Ghosh}}.\;$K$-frame is also presented to reconstruct elements from the range of a bounded linear operator \,$K$\, in a separable Hilbert space and it is a generalization of the ordinary frames.\,Infact, many properties of ordinary frames may not holds for such generalization of frames.\,Like \,$K$-frame, another generalization of frame is \,$g$-fusion frame and it has been studied by several authors {\cite{P, Sadri, Ahmadi}}.\,S.\,Rabinson \cite{S} presented the basic concepts of tensor product of Hilbert spaces.\,The tensor product of Hilbert spaces \,$X$\, and \,$Y$\, is a certain linear space of operators which was represented by Folland in \cite{Folland}, Kadison and Ringrose in \cite{Kadison}.\,Generalized fusion frame in tensor product of Hilbert spaces was studied by P.\,Ghosh and T.\,K.\,Samanta \cite{K}.

In 1970, Diminnie et. al. \cite{Diminnie} introduced the concept of \,$2$-inner product space. Atomic system in \,$2$-inner product space is studied by D. Bahram and J. Mohammad \cite{Janfada}.\;A generalization of a \,$2$-inner product space for \,$n \,\geq\, 2$\, was developed by A.\,Misiak \cite{Misiak} in 1989.

In this paper, we give a notion of a family of local atoms in \,$n$-Hilbert space.\,Since tensor product of \,$n$-Hilbert spaces becomes a \,$n$-Hilbert space, we like to study \,$K$-frame in this \,$n$-Hilbert space.\,We give a necessary and sufficient condition for being \,$K$-frames in \,$n$-Hilbert spaces is that of being in their tensor products.\,Atomic system in tensor product of \,$n$-Hilbert spaces is discussed.\,Finally, we are going to establish a relationship between atomic systems in \,$n$-Hilbert spaces and their tensor products. 

Throughout this paper,\;$X$\, will denote separable Hilbert spaces with the inner product \;$\left <\,\cdot \,,\, \cdot\, \right>_{1}$\, and \,$\mathbb{K}$\, denote the field of real or complex numbers.\;$l^{\,2}\,(\,\mathbb{N}\,)$\, and $l^{\,2}\,(\,\mathbb{N} \,\times\, \mathbb{N}\,)$\, denote the spaces of square summable scalar-valued sequences with  index sets \,$\mathbb{N}$\, and \,$\mathbb{N} \,\times\, \mathbb{N}$, respectively.\;$\mathcal{B}\,(\,X\,)$\, denote the space of all bounded linear operators on \,$X$.

\section{Preliminaries}
\smallskip\hspace{.6 cm}

\begin{definition}\cite{L}
Let \,$K \,\in\, \mathcal{B}\,(\,X\,)$.\;A sequence \,$\left\{\,f_{\,i}\,\right\}_{i \,=\, 1}^{\infty} \,\subseteq\, X$\, is called a \,$K$-frame for \,$X$\, if there exist positive constants \,$A,\, B$\, such that
\begin{equation}\label{eq0.01}
 A\; \|\,K^{\,\ast}\,f\,\|_{1}^{\,2} \,\leq\, \;\sum\limits_{i \,=\, 1}^{\infty}\;  \left|\ \left <\,f\,,\, f_{\,i} \, \right >_{1}\,\right|^{\,2} \,\leq\, B \;\|\,f\,\|_{1}^{\,2}\; \;\forall\; f \,\in\, X
\end{equation}
The constants \,$A,\,B$\, are called frame bounds.\,If \,$\left\{\,f_{\,i}\,\right\}_{i \,=\, 1}^{\infty}$\; satisfies only 
\[\left|\ \left <\,f\,,\, f_{\,i} \, \right >_{1}\,\right|^{\,2} \,\leq\, B \;\|\,f\,\|_{1}^{\,2}\; \;\forall\; f \,\in\, X\]
then it is called a Bessel sequence with bound \,$B$.
\end{definition}

\begin{definition}\cite{L}
Let \,$K \,\in\, \mathcal{B}\,(\,X\,)$\, and \,$\left\{\,f_{\,i}\,\right\}_{i \,=\, 1}^{\infty}$\, be a sequence in \,$X$.\,Then \,$\left\{\,f_{\,i}\,\right\}_{i \,=\, 1}^{\infty}$\, is  said to be an atomic system for \,$K$\, if the following statements hold:
\begin{itemize}
\item[(I)] the series \,$\sum\limits_{i \,=\, 1}^{\,\infty}\,c_{\,i}\,f_{\,i}$\, converges for all \,$\{\,c_{\,i}\,\} \,\in\, l^{\,2}\,(\,\mathbb{N}\,)$;
\item[(II)] for every \,$x \,\in\, X$, there exists \,$a_{\,x} \,=\, \{\,a_{\,i}\,\} \,\in\, l^{\,2}\,(\,\mathbb{N}\,)$\, such that \,$\left\|\,a_{\,x}\,\right\|_{l^{\,2}} \,\leq\, C\,\|\,x\,\|_{1}$\, and \,$K\,(\,x\,) \,=\, \sum\limits_{i}\,a_{\,i}\,f_{\,i}$, for some \,$C \,>\, 0$.
\end{itemize} 
\end{definition}

\begin{definition}\cite{Upender}\label{def0.001}
Let \,$\left(\,Y,\, \left<\,\cdot,\, \cdot\,\right>_{2}\,\right)$\, be a Hilbert space.\,Then the tensor product of \,$X$\, and \,$Y$\,  is denoted by \,$X \,\otimes\, Y$\, and it is defined to be an inner product space associated with the inner product   
\[\left<\,f \,\otimes\, g\,,\, f^{\,\prime} \,\otimes\, g^{\,\prime}\,\right> \,=\, \left<\,f\,,\, f^{\,\prime}\,\right>_{1}\;\left<\,g\,,\, g^{\,\prime}\,\right>_{2}\; \;\forall\; f,\, f^{\,\prime} \,\in\, X\; \;\&\; \;g,\, g^{\,\prime} \,\in\, Y.\]
The norm on \,$X \,\otimes\, Y$\, is given by 
\[\left\|\,f \,\otimes\, g\,\right\| \,=\, \|\,f\,\|_{\,1}\;\|\,g\,\|_{\,2}\; \;\forall\; f \,\in\, X\; \;\&\; \,g \,\in\, Y.\]
The space \,$X \,\otimes\, Y$\, is a Hilbert space with respect to the above inner product.     
\end{definition} 

For \,$Q \,\in\, \mathcal{B}\,(\,H\,)$\, and \,$T \,\in\, \mathcal{B}\,(\,K\,)$, the tensor product of operators \,$Q$\, and \,$T$\, is denoted by \,$Q \,\otimes\, T$\, and defined as 
\[\left(\,Q \,\otimes\, T\,\right)\,A \,=\, Q\,A\,T^{\,\ast}\; \;\forall\; \;A \,\in\, H \,\otimes\, K.\]

\begin{theorem}\cite{Folland,Li}\label{th1.1}
Suppose \,$Q,\, Q^{\prime} \,\in\, \mathcal{B}\,(\,X\,)$\, and \,$T,\, T^{\prime} \,\in\, \mathcal{B}\,(\,Y\,)$, then \begin{itemize}
\item[(I)]\hspace{.2cm} \,$Q \,\otimes\, T \,\in\, \mathcal{B}\,(\,X \,\otimes\, Y\,)$\, and \,$\left\|\,Q \,\otimes\, T\,\right\| \,=\, \|\,Q\,\|\; \|\,T\,\|$.
\item[(II)]\hspace{.2cm} \,$\left(\,Q \,\otimes\, T\,\right)\,(\,f \,\otimes\, g\,) \,=\, Q\,f \,\otimes\, T\,g$\, for all \,$f \,\in\, H,\, g \,\in\, K$.
\item[(III)]\hspace{.2cm} $\left(\,Q \,\otimes\, T\,\right)\,\left(\,Q^{\,\prime} \,\otimes\, T^{\,\prime}\,\right) \,=\, (\,Q\,Q^{\,\prime}\,) \,\otimes\, (\,T\,T^{\,\prime}\,)$. 
\item[(IV)]\hspace{.2cm} \,$Q \,\otimes\, T$\, is invertible if and only if \,$Q$\, and \,$T$\, are invertible, in which case \,$\left(\,Q \,\otimes\, T\,\right)^{\,-\, 1} \,=\, \left(\,Q^{\,-\, 1} \,\otimes\, T^{\,-\, 1}\,\right)$.
\item[(V)]\hspace{.2cm} \,$\left(\,Q \,\otimes\, T\,\right)^{\,\ast} \,=\, \left(\,Q^{\,\ast} \,\otimes\, T^{\,\ast}\,\right)$.  
\item[(VI)]\hspace{.2cm} Let \,$f,\, f^{\,\prime} \,\in\, X \,\setminus\, \{\,0\}$\, and \,$g,\, g^{\,\prime} \,\in\, Y \,\setminus\, \{\,0\}$\,.\;If \,$f \,\otimes\, g \,=\, f^{\,\prime} \,\otimes\, g^{\,\prime}$, then there exist constants \,$A$\, and \,$B$\, with \,$A\,B \,=\, 1$\, such that \,$f \,=\, A\,f^{\,\prime}$\, and \,$g \,=\, B\,g^{\,\prime}$.  
\end{itemize}
\end{theorem}

\begin{definition}\cite{Mashadi}
A \,$n$-norm on a linear space \,$H$\, over the field \,$\mathbb{K}$\, is a function
\[\left(\,x_{\,1},\, x_{\,2},\, \cdots,\, x_{\,n}\,\right) \,\longmapsto\, \left\|\,x_{\,1},\, x_{\,2},\, \cdots,\, x_{\,n}\,\right\|,\; x_{\,1},\, x_{\,2},\, \cdots,\, x_{\,n} \,\in\, H\]from \,$H^{\,n}$\, to the set \,$\mathbb{R}$\, of all real numbers such that for every \,$x_{\,1},\, x_{\,2},\, \cdots,\, x_{\,n} \,\in\, H$,
\begin{itemize}
\item[(I)]\;\; $\left\|\,x_{\,1},\, x_{\,2},\, \cdots,\, x_{\,n}\,\right\| \,=\, 0$\; if and only if \,$x_{\,1},\, \cdots,\, x_{\,n}$\; are linearly dependent,
\item[(II)]\;\;\; $\left\|\,x_{\,1},\, x_{\,2},\, \cdots,\, x_{\,n}\,\right\|$\, is invariant under permutations of \,$x_{\,1},\, x_{\,2},\, \cdots,\, x_{\,n}$,
\item[(III)]\;\;\; $\left\|\,\alpha\,x_{\,1},\, x_{\,2},\, \cdots,\, x_{\,n}\,\right\| \,=\, |\,\alpha\,|\, \left\|\,x_{\,1},\, x_{\,2},\, \cdots,\, x_{\,n}\,\right\|$, for \,$\alpha \,\in\, \mathbb{K}$,
\item[(IV)]\;\; $\left\|\,x \,+\, y,\, x_{\,2},\, \cdots,\, x_{\,n}\,\right\| \,\leq\, \left\|\,x,\, x_{\,2},\, \cdots,\, x_{\,n}\,\right\| \,+\,  \left\|\,y,\, x_{\,2},\, \cdots,\, x_{\,n}\,\right\|$.
\end{itemize}
A linear space \,$H$, together with a n-norm \,$\left\|\,\cdot,\, \cdots,\, \cdot \,\right\|$, is called a linear n-normed space. 
\end{definition}

\begin{definition}\cite{Misiak}
Let \,$n \,\in\, \mathbb{N}$\; and \,$H$\, be a linear space of dimension greater than or equal to \,$n$\; over the field \,$\mathbb{K}$.\;An n-inner product on \,$H$\, is a map 
\[\left(\,x,\, y,\, x_{\,2},\, \cdots,\, x_{\,n}\,\right) \,\longmapsto\, \left<\,x,\, y \,|\, x_{\,2},\, \cdots,\, x_{\,n} \,\right>,\; x,\, y,\, x_{\,2},\, \cdots,\, x_{\,n} \,\in\, H\]from \,$H^{n \,+\, 1}$\, to the set \,$\mathbb{K}$\, such that for every \,$x,\, y,\, x_{\,1},\, x_{\,2},\, \cdots,\, x_{\,n} \,\in\, H$,
\begin{itemize}
\item[(I)]\;\; $\left<\,x_{\,1},\, x_{\,1} \,|\, x_{\,2},\, \cdots,\, x_{\,n} \,\right> \,\geq\,  0$\, and \,$\left<\,x_{\,1},\, x_{\,1} \,|\, x_{\,2},\, \cdots,\, x_{\,n} \,\right> \,=\,  0$\, if and only if \,$x_{\,1},\, x_{\,2},\, \cdots,\, x_{\,n}$\, are linearly dependent,
\item[(II)]\;\; $\left<\,x,\, y \,|\, x_{\,2},\, \cdots,\, x_{\,n} \,\right> \,=\, \left<\,x,\, y \,|\, x_{\,i_{\,2}},\, \cdots,\, x_{\,i_{\,n}}\,\right> $\, for every permutations \,$\left(\, i_{\,2},\, \cdots,\, i_{\,n} \,\right)$\, of \,$\left(\, 2,\, \cdots,\, n \,\right)$,
\item[(III)]\;\; $\left<\,x,\, y \,|\, x_{\,2},\, \cdots,\, x_{\,n} \,\right> \,=\, \overline{\left<\,y,\, x \,|\, x_{\,2},\, \cdots,\, x_{\,n} \,\right> }$,
\item[(IV)]\;\; $\left<\,\alpha\,x,\, y \,|\, x_{\,2},\, \cdots,\, x_{\,n} \,\right> \,=\, \alpha \,\left<\,x,\, y \,|\, x_{\,2},\, \cdots,\, x_{\,n} \,\right> $, for \,$\alpha \,\in\, \mathbb{K}$,
\item[(V)]\;\; $\left<\,x \,+\, y,\, z \,|\, x_{\,2},\, \cdots,\, x_{\,n} \,\right> \,=\, \left<\,x,\, z \,|\, x_{\,2},\, \cdots,\, x_{\,n}\,\right> \,+\,  \left<\,y,\, z \,|\, x_{\,2},\, \cdots,\, x_{\,n} \,\right>$.
\end{itemize}
A linear space \,$H$\, together with an n-inner product \,$\left<\,\cdot,\, \cdot \,|\, \cdot,\, \cdots,\, \cdot\,\right>$\, is called an n-inner product space.
\end{definition}

\begin{theorem}(Schwarz inequality)\cite{Misiak}\label{thn1}
Let \,$H$\, be a \,$n$-inner product space.\,Then 
\[\left|\,\left<\,x,\, y \,|\, x_{\,2},\, \cdots,\, x_{\,n}\,\right>\,\right| \,\leq\, \left\|\,x,\, x_{\,2},\, \cdots,\, x_{\,n}\,\right\|\, \left\|\,y,\, x_{\,2},\, \cdots,\, x_{\,n}\,\right\|\]
hold for all \,$x,\, y,\, x_{\,2},\, \cdots,\, x_{\,n} \,\in\, H$.
\end{theorem}

\begin{theorem}\cite{Misiak}
Let \,$H$\, be a \,$n$-inner product space.\,Then
\[\left \|\,x_{\,1},\, x_{\,2},\, \cdots,\, x_{\,n}\,\right\| \,=\, \sqrt{\left <\,x_{\,1},\, x_{\,1} \,|\, x_{\,2},\, \cdots,\, x_{\,n}\,\right>}\] defines a n-norm for which
\begin{align*}
&\|\,x \,+\, y,\, x_{\,2},\, \cdots,\, x_{\,n}\,\|^{\,2} \,+\, \|\,x \,-\, y,\, x_{\,2},\, \cdots,\, x_{\,n}\,\|^{\,2}\\
& \,=\, 2\, \left(\,\|\,x,\, x_{\,2},\, \cdots,\, x_{\,n}\,\|^{\,2} \,+\, \|\,y,\, x_{\,2},\, \cdots,\, x_{\,n}\,\|^{\,2} \,\right)
\end{align*} 
hold for all \,$x,\, y,\, x_{\,1},\, x_{\,2},\, \cdots,\, x_{\,n} \,\in\, H$.
\end{theorem}

\begin{definition}\cite{Mashadi}
Let \,$\left(\,H,\, \left\|\,\cdot,\, \cdots,\, \cdot \,\right\|\,\right)$\; be a linear n-normed space.\;A sequence \,$\{\,x_{\,k}\,\}$\; in \,$H$\, is said to convergent if there exists an \,$x \,\in\, H$\, such that 
\[\lim\limits_{k \to \infty}\,\left\|\,x_{\,k} \,-\, x \,,\, x_{\,2} \,,\, \cdots \,,\, x_{\,n} \,\right\| \,=\, 0\]
for every \,$ x_{\,2},\, \cdots,\, x_{\,n} \,\in\, H$\, and it is called a Cauchy sequence if 
\[\lim\limits_{l \,,\, k \,\to\, \infty}\,\left \|\,x_{l} \,-\, x_{\,k} \,,\, x_{\,2} \,,\, \cdots \,,\, x_{\,n}\,\right\| \,=\, 0\]
for every \,$ x_{\,2},\, \cdots,\, x_{\,n} \,\in\, H$.\;The space \,$H$\, is said to be complete if every Cauchy sequence in this space is convergent in \,$H$.\;A n-inner product space is called n-Hilbert space if it is complete with respect to its induce norm.
\end{definition}


\section{Atomic system in $n$-Hilbert space}

\smallskip\hspace{.6 cm} In this section, concept of a family of local atoms associated to \,$\left(\,a_{\,2},\, \cdots,\, a_{\,n}\,\right)$\, is discussed.\;Next, we are going to generalize this concept and then define \,$K$-frame associated to \,$\left(\,a_{\,2},\, \cdots,\, a_{\,n}\,\right)$\, for \,$H$, for a given bounded linear operator \,$K$. \\

Let \,$a_{\,2},\, a_{\,3},\, \cdots,\, a_{\,n}$\, be the fixed elements in \,$H$\, and \,$L_{F}$\, denote the linear subspace of \,$H$\, spanned by the non-empty finite set \,$F \,=\, \left\{\,\,a_{\,2} \,,\, a_{\,3} \,,\, \cdots \,,\, a_{\,n}\,\right\}$. Then the quotient space \,$H \,/\, L_{F}$\, is a normed linear space with respect to the norm, \,$\left\|\,x \,+\, L_{F}\,\right\|_{F} \,=\, \left\|\,x \,,\, a_{\,2} \,,\,  \cdots \,,\, a_{\,n}\,\right\|$\, for every \,$x \,\in\, H$.\;Let \,$M_{F}$\, be the algebraic complement of \,$L_{F}$, then \,$H \,=\, L_{F} \,\oplus\, M_{F}$.\;Define   
\[\left<\,x \,,\, y\,\right>_{F} \,=\, \left<\,x \,,\, y \;|\; a_{\,2} \,,\,  \cdots \,,\, a_{\,n}\,\right>\; \;\text{on}\; \;H.\]
Then \,$\left<\,\cdot \,,\, \cdot\,\right>_{F}$\, is a semi-inner product on \,$H$\, and this semi-inner product induces an inner product on the quotient space \,$H \,/\, L_{F}$\; which is given by
\[\left<\,x \,+\, L_{F} \,,\, y \,+\, L_{F}\,\right>_{F} \,=\, \left<\,x \,,\, y\,\right>_{F} \,=\, \left<\,x \,,\, y \,|\, a_{\,2} \,,\,  \cdots \,,\, a_{\,n} \,\right>\;\; \;\forall \;\; x,\, y \,\in\, H.\]
By identifying \,$H \,/\, L_{F}$\; with \,$M_{F}$\; in an obvious way, we obtain an inner product on \,$M_{F}$.\;Then \,$M_{F}$\, is a normed space with respect to the norm \,$\|\,\cdot\,\|_{F}$\, defined by \,$\|\,x\,\|_{F} \;=\; \sqrt{\left<\,x \,,\, x \,\right>_{F}}\; \;\forall\, x \,\in\, M_{F}$.\;Let \,$H_{F}$\, be the completion of the inner product space \,$M_{F}$.

\begin{definition}\cite{Prasenjit}\label{def0.1}
Let \,$H$\, be a n-Hilbert space.\;A sequence \,$\left\{\,f_{\,i}\,\right\}^{\,\infty}_{\,i \,=\, 1}$\, in \,$H$\, is said to be a frame associated to \,$\left(\,a_{\,2},\, \cdots,\, a_{\,n}\,\right)$\, if there exists constant \,$0 \,<\, A \,\leq\, B \,<\, \infty$\,  such that  
\[ A \, \left\|\,f \,,\, a_{\,2} \,,\, \cdots \,,\, a_{\,n} \,\right\|^{\,2} \,\leq\, \sum\limits^{\infty}_{i \,=\, 1}\,\left|\,\left<\,f \,,\, f_{\,i} \,|\, a_{\,2} \,,\, \cdots \,,\, a_{\,n}\,\right>\,\right|^{\,2} \,\leq\, B\, \left\|\,f \,,\, a_{\,2} \,,\, \cdots \,,\, a_{\,n}\,\right\|^{\,2}\]for all \,$f \,\in\, H$.\,The constants \,$A,\,B$\, are called frame bounds.\,If \,$\left\{\,f_{\,i}\,\right\}^{\,\infty}_{\,i \,=\, 1}$\; satisfies 
\[\sum\limits^{\infty}_{i \,=\, 1}\,\left|\,\left<\,f \,,\, f_{\,i} \,|\, a_{\,2} \,,\, \cdots \,,\, a_{\,n}\,\right>\,\right|^{\,2} \,\leq\, B\; \left\|\,f \,,\, a_{\,2} \,,\, \cdots \,,\, a_{\,n}\,\right\|^{\,2}\; \;\forall\; f \,\in\, H\] is called a Bessel sequence associated to \,$\left(\,a_{\,2},\, \cdots,\, a_{\,n}\,\right)$\, in \,$H$\, with bound \,$B$.
\end{definition}

\begin{theorem}\label{th2}
Let \,$H$\, be a n-Hilbert space.\;Then \,$\left\{\,f_{\,i}\,\right\}^{\,\infty}_{\,i \,=\, 1} \,\subseteq\, H$\; is a frame associated to \,$\left(\,a_{\,2},\, \cdots,\, a_{\,n}\,\right)$\; with bounds \,$A \;\;\&\;\; B$\; if and only if it is a frame for the Hilbert space \,$H_{F}$\; with bounds \,$A \;\;\&\;\; B$.
\end{theorem}

\begin{proof}
This theorem is an extension of the Theorem (3.2) of \cite{Sadeghi} and proof of this theorem directly follows from the Theorem (3.2) of \cite{Sadeghi}.
\end{proof}
For more details on frames in \,$n$-Hilbert spaces and their tensor products one can go through the papers \cite{Prasenjit,G}.
  
\begin{definition}
Let \,$\left(\,H \,,\, \|\,\cdot \,,\, \cdots \,,\, \cdot\,\| \,\right)$\; be a linear n-normed space and \,$a_{\,2},\, \cdots,\, a_{\,n}$\, be fixed elements in \,$H$.\;Let \,$W$\, be a subspace of \,$H$\, and \,$\left<\,a_{\,i}\,\right>$\, denote the subspaces of \,$H$\, generated by \,$a_{\,i}$, for \,$i \,=\, 2,\, 3,\, \cdots,\,n $.\;Then a map \,$T \,:\, W \,\times\,\left<\,a_{\,2}\,\right> \,\times\, \cdots \,\times\, \left<\,a_{\,n}\,\right> \,\to\, \mathbb{K}$\; is called a b-linear functional on \,$W \,\times\, \left<\,a_{\,2}\,\right> \,\times\, \cdots \,\times\, \left<\,a_{\,n}\,\right>$, if for every \,$x,\, y \,\in\, W$\, and \,$k \,\in\, \mathbb{K}$, the following conditions hold:
\begin{itemize}
\item[(I)]\hspace{.2cm} $T\,(\,x \,+\, y \,,\, a_{\,2}  \,,\, \cdots \,,\, a_{\,n}\,) \,=\, T\,(\,x  \,,\, a_{\,2} \,,\, \cdots \,,\, a_{\,n}\,) \,+\, T\,(\,y  \,,\, a_{\,2} \,,\, \cdots \,,\, a_{\,n}\,)$
\item[(II)]\hspace{.2cm} $T\,(\,k\,x  \,,\, a_{\,2} \,,\, \cdots \,,\, a_{\,n}\,) \,=\, k\; T\,(\,x  \,,\, a_{\,2} \,,\, \cdots \,,\, a_{\,n}\,)$. 
\end{itemize}
A b-linear functional is said to be bounded if \,$\exists$\, a real number \,$M \,>\, 0$\; such that
\[\left|\,T\,(\,x  \,,\, a_{\,2} \,,\, \cdots \,,\, a_{\,n}\,)\,\right| \,\leq\, M\; \left\|\,x  \,,\, a_{\,2} \,,\, \cdots \,,\, a_{\,n}\,\right\|\; \;\forall\; x \,\in\, W.\]
\end{definition}

Some properties of bounded\;$b$-linear functional defined on \,$H \,\times\, \left<\,a_{\,2}\,\right> \,\times\, \cdots \,\times\, \left<\,a_{\,n}\,\right>$\, have been discussed in \cite{T}.

\begin{definition}
Let \,$\left\{\,f_{\,i}\,\right\}_{i \,=\, 1}^{\,\infty}$\, be a Bessel sequence associated to \,$\left(\,a_{\,2},\,\cdots,\, a_{\,n}\,\right)$\, in \,$H$\, and \,$Y$\, be a closed subspace of \,$H$.\;Then \,$\left\{\,f_{\,i}\,\right\}_{i \,=\, 1}^{\,\infty}$\, is said to be a family of local atoms associated to \,$\left(\,a_{\,2},\, \cdots,\, a_{\,n}\,\right)$\, for \,$Y$\, if there exists a sequence of bounded b-linear functionals \,$\left\{\,T_{\,i}\,\right\}_{i \,=\, 1}^{\,\infty}$\, defined on \,$H \,\times\, \left<\,a_{\,2}\,\right> \,\times\, \,\cdots\, \times\, \left<\,a_{\,n}\,\right>$\, such that
\begin{itemize}
\item[(I)] \,$\sum\limits_{i \,=\, 1}^{\,\infty}\,\left|\,T_{i}\,\left(\,f,\, a_{\,2},\, \cdots,\, a_{\,n}\,\right)\,\right|^{\,2} \,\leq\, C\,\left\|\,f \,,\, a_{\,2} \,,\, \cdots \,,\, a_{\,n}\,\right\|^{\,2}$, for some \,$C \,>\, 0$.
\item[(II)] \,$f \,=\, \sum\limits_{i \,=\, 1}^{\,\infty}\,T_{i}\,\left(\,f,\, a_{\,2},\, \cdots,\, a_{\,n}\,\right)\,f_{\,i}$, for all \,$f \,\in\, Y$.
\end{itemize}  
\end{definition}

\begin{theorem}
Let \,$\left\{\,f_{\,i}\,\right\}_{i \,=\, 1}^{\,\infty}$\, be a family of local atoms associated to \,$\left(\,a_{\,2},\, \cdots,\, a_{\,n}\,\right)$\, for \,$Y$, where \,$Y$\, be a closed subspace of \,$H$.\;Then the family \,$\left\{\,f_{\,i}\,\right\}_{i \,=\, 1}^{\,\infty}$\, is a frame associated to \,$\left(\,a_{\,2},\, \cdots,\, a_{\,n}\,\right)$\, for \,$Y$. 
\end{theorem}

\begin{proof}
Since \,$\left\{\,f_{\,i}\,\right\}_{i \,=\, 1}^{\,\infty}$\, is a family of local atoms associated to \,$\left(\,a_{\,2},\, \cdots,\, a_{\,n}\,\right)$\, for \,$Y$, there exists a sequence of bounded \,$b$-linear functionals \,$\left\{\,T_{\,i}\,\right\}_{i \,=\, 1}^{\,\infty}$\, such that 
\[\sum\limits_{i \,=\, 1}^{\,\infty}\,\left|\,T_{i}\,\left(\,f,\, a_{\,2},\, \cdots,\, a_{\,n}\,\right)\,\right|^{\,2} \,\leq\, C\,\left\|\,f \,,\, a_{\,2} \,,\, \cdots \,,\, a_{\,n}\,\right\|^{\,2}, f \,\in\, Y, \,\text{for some \,$C \,>\, 0$}\]
Now, for each \,$f \,\in\, Y$,
\[\left\|\,f \,,\, a_{\,2} \,,\, \cdots \,,\, a_{\,n}\,\right\|^{\,4} \,=\, \left(\,\left<\,f \,,\, f \,|\, a_{\,2} \,,\, \cdots \,,\, a_{\,n}\,\right>\,\right)^{\,2}\hspace{1.5cm}\]
\[=\, \left(\,\left<\,f \,,\, \sum\limits_{i \,=\, 1}^{\,\infty}\,T_{i}\,\left(\,f,\, a_{\,2},\, \cdots,\, a_{\,n}\,\right)\,f_{\,i} \,|\, a_{\,2} \,,\, \cdots \,,\, a_{\,n}\,\right>\,\right)^{\,2}\]
\[=\, \left(\,\sum\limits_{i \,=\, 1}^{\,\infty}\,\overline{\,T_{i}\,\left(\,f,\, a_{\,2},\, \cdots,\, a_{\,n}\,\right)\,}\,\left<\,f \,,\, f_{\,i} \,|\, a_{\,2} \,,\, \cdots \,,\, a_{\,n}\,\right>\,\right)^{\,2}\hspace{.1cm}\] 
\[\hspace{.5cm}\leq\, \sum\limits_{i \,=\, 1}^{\,\infty}\,\left|\,T_{i}\,\left(\,f,\, a_{\,2},\, \cdots,\, a_{\,n}\,\right)\,\right|^{\,2}\,\sum\limits_{i \,=\, 1}^{\,\infty}\, \left|\,\left <\,f \,,\, f_{\,i} \,|\, a_{\,2} \,,\, \cdots \,,\, a_{\,n}\,\right >\,\right|^{\,2}\]
\[\leq\, C\,\left\|\,f \,,\, a_{\,2} \,,\, \cdots \,,\, a_{\,n}\,\right\|^{\,2}\,\sum\limits_{i \,=\, 1}^{\,\infty}\, \left|\,\left <\,f \,,\,  f_{\,i} \,|\, a_{\,2} \,,\, \cdots \,,\, a_{\,n}\,\right >\,\right|^{\,2}\hspace{.3cm}\]
\[\hspace{.5cm}\Rightarrow\, \dfrac{1}{C}\,\left\|\,f \,,\, a_{\,2} \,,\, \cdots \,,\, a_{\,n}\,\right\|^{\,2} \,\leq\, \sum\limits_{i \,=\, 1}^{\,\infty}\, \left|\,\left <\,f \,,\, f_{\,i} \,|\, a_{\,2} \,,\, \cdots \,,\, a_{\,n}\,\right >\,\right|^{\,2}.\]
Also, \,$\left\{\,f_{\,i}\,\right\}_{i \,=\, 1}^{\,\infty}$\, is a Bessel sequence associated to \,$\left(\,a_{\,2},\, \cdots,\, a_{\,n}\,\right)$\, in \,$Y$.\,Hence, \,$\left\{\,f_{\,i}\,\right\}_{i \,=\, 1}^{\,\infty}$\, is a frame associated to \,$\left(\,a_{\,2},\, \cdots,\, a_{\,n}\,\right)$\, for \,$Y$.    
\end{proof}

\begin{theorem}
Let \,$\left\{\,f_{\,i}\,\right\}_{i \,=\, 1}^{\,\infty}$\, be a Bessel sequence associated to \,$\left(\,a_{\,2},\,\cdots,\, a_{\,n}\,\right)$\, in \,$H$\, and \,$Y$\, be a closed subspace of \,$H$.\;If there exists a Bessel sequence associated to \,$\left(\,a_{\,2},\,\cdots,\, a_{\,n}\,\right)$\, in \,$H$, say \,$\left\{\,g_{\,i}\,\right\}_{i \,=\, 1}^{\,\infty}$\, such that
\begin{equation}\label{eqm1} 
P_{Y}\,(\,f\,) \,=\, \sum\limits_{i \,=\, 1}^{\,\infty}\,\left <\,f \,,\,  g_{\,i} \,|\, a_{\,2} \,,\, \cdots \,,\, a_{\,n}\,\right >\,f_{\,i}, \;\;\text{for all \,$f \,\in\, H_{F}$},
\end{equation}
where \,$P_{Y}$\, is the orthogonal projection onto \,$Y$, then \,$\left\{\,f_{\,i}\,\right\}_{i \,=\, 1}^{\,\infty}$\, is a family of local atoms associated to \,$\left(\,a_{\,2},\, \cdots,\, a_{\,n}\,\right)$\, for \,$Y$. 
\end{theorem}

\begin{proof}
Let us take \,$f \,\in\, Y$, then by (\ref{eqm1}), we can write
\[f \,=\, P_{Y}\,(\,f\,) \,=\, \sum\limits_{i \,=\, 1}^{\,\infty}\,\left<\,f \,,\,  g_{\,i} \,|\, a_{\,2} \,,\, \cdots \,,\, a_{\,n}\,\right>\,f_{\,i}.\]
\[\text{Now, we define}\hspace{.5cm} T_{i}\,\left(\,f,\, a_{\,2},\, \cdots,\, a_{\,n}\,\right) \,=\, \left<\,f \,,\,  g_{\,i} \,|\, a_{\,2} \,,\, \cdots \,,\, a_{\,n}\,\right>\; \;\forall\; f \,\in\, Y.\hspace{3cm}\] 
\[\text{Then}\hspace{.5cm} f \,=\, \sum\limits_{i \,=\, 1}^{\,\infty}\,T_{i}\,\left(\,f,\, a_{\,2},\, \cdots,\, a_{\,n}\,\right)\,f_{\,i}\; \;\forall\; f \,\in\, Y.\hspace{7cm}\]Also, for any \,$i$, we have
\[\left|\,T_{i}\,\left(\,f,\, a_{\,2},\, \cdots,\, a_{\,n}\,\right)\,\right| \,=\, \left|\,\left<\,f \,,\,  g_{\,i} \,|\, a_{\,2} \,,\, \cdots \,,\, a_{\,n}\,\right>\,\right|\hspace{3cm}\]
\[\leq\,\left\|\,f,\, a_{\,2},\, \cdots,\, a_{\,n}\,\right\|\,\left\|\,g_{\,i},\, a_{\,2},\, \cdots,\, a_{\,n}\,\right\|\; \;[\;\text{by Theorem (\ref{thn1})}\;]\hspace{1cm}\]
\[\leq\, M\,\left\|\,f,\, a_{\,2},\, \cdots,\, a_{\,n}\,\right\|,\; \left[\;\text{where \,$M \,=\, \sup\limits_{i}\,\left\|\,g_{\,i},\, a_{\,2},\, \cdots,\, a_{\,n}\,\right\|$}\right].\hspace{.3cm}\]
This verifies that each \,$T_{\,i}$\, are bounded \,$b$-linear functionals defined on \,$Y \,\times\, \left<\,a_{\,2}\,\right> \,\times\, \,\cdots\, \times\, \left<\,a_{\,n}\,\right>$.\,On the other hand, we get 
\[\sum\limits_{i \,=\, 1}^{\,\infty}\left|\,T_{i}\left(\,f,\, a_{\,2},\, \cdots,\, a_{\,n}\,\right)\,\right|^{\,2} = \sum\limits_{i \,=\, 1}^{\,\infty}\left|\,\left<\,f,\,  g_{\,i} \,|\, a_{\,2},\, \cdots,\, a_{\,n}\,\right>\,\right|^{\,2}\leq\, B\left\|\,f,\, a_{\,2},\, \cdots,\, a_{\,n}\,\right\|^{\,2}\]
\[\left[\,\text{since \,$\left\{\,g_{\,i}\,\right\}_{i \,=\, 1}^{\,\infty}$\, is a Bessel sequence associated to \,$\left(\,a_{\,2},\,\cdots,\, a_{\,n}\,\right)$}\,\right].\]
This completes the proof. 
\end{proof}

Now, we are going to generalize the concept of a family of local atoms associated to \,$\left(\,a_{\,2},\,\cdots,\, a_{\,n}\,\right)$.

\begin{definition}\label{defn1}
Let \,$K$\, be a bounded linear operator on \,$H_{F}$\, and \,$\left\{\,f_{\,i}\,\right\}_{i \,=\, 1}^{\,\infty}$\, be a sequence of vectors in \,$H$.\;Then \,$\left\{\,f_{\,i}\,\right\}_{i \,=\, 1}^{\,\infty}$\, is said to be an atomic system associated to \,$\left(\,a_{\,2},\, \cdots,\, a_{\,n}\,\right)$\, for \,$K$\, in \,$H$\, if 
\begin{itemize}
\item[(I)]\,$\left\{\,f_{\,i}\,\right\}_{i \,=\, 1}^{\,\infty}$\, is a Bessel sequence associated to \,$\left(\,a_{\,2},\,\cdots,\, a_{\,n}\,\right)$\, in \,$H$.
\item[(II)] For any \,$f \,\in\, H_{F}$, there exists \,$\{\,c_{\,i}\,\}_{i \,=\, 1}^{\,\infty} \,\in\, l^{\,2}\,(\,\mathbb{N}\,)$\, such that \,$K\,(\,f\,) =\, \sum\limits_{i \,=\, 1}^{\,\infty}\,c_{\,i}\,f_{\,i}$, where \,$\left\|\,\{\,c_{\,i}\,\}_{i \,=\, 1}^{\,\infty}\,\right\|_{l^{2}} \,\leq\, C\, \left\|\,f \,,\, a_{\,2} \,,\, \cdots \,,\, a_{\,n} \,\right\|$\, and \,$C \,>\, 0$. 
\end{itemize} 
\end{definition}

\begin{definition}
Let \,$K$\, be a bounded linear operator on \,$H_{F}$.\,Then a sequence \,$\{\,f_{\,i}\,\}_{i \,=\, 1}^{\infty} \,\subseteq\, H$\, is said to be a \,$K$-frame associated to \,$\left(\,a_{\,2},\, \cdots,\, a_{\,n}\,\right)$\, for \,$H$\, if there exist constants \,$A,\, B \,>\, 0$\, such that for each \,$f \,\in\, H_{F}$, 
\[A\,\left \|\,K^{\,\ast}\,f,\, a_{\,2},\, \cdots,\, a_{\,n}\,\right \|^{\,2} \,\leq\, \sum\limits_{i \,=\, 1}^{\infty}\, \left|\,\left <\,f,\,  f_{\,i} \,|\, a_{\,2},\, \cdots,\, a_{\,n}\,\right >\,\right|^{\,2} \,\leq\, B\,\left \|\,f,\, a_{\,2},\, \cdots,\, a_{\,n}\,\right \|^{\,2}.\]
\end{definition}

\begin{theorem}\label{th2.1}
Let \,$\{\,f_{\,i}\,\}_{i \,=\, 1}^{\infty}$\, be a \,$K$-frame associated to \,$\left(\,a_{\,2},\, \cdots,\, a_{\,n}\,\right)$\, for \,$H$.\;Then there exists a Bessel sequence \,$\{\,g_{\,i}\,\}_{i \,=\, 1}^{\infty}$\, associated to \,$\left(\,a_{\,2},\, \cdots,\, a_{\,n}\,\right)$\, such that
\[K^{\,\ast}\,f \,=\, \sum\limits_{i \,=\, 1}^{\infty}\,\left <\,f \,,\,  f_{\,i} \,|\, a_{\,2} \,,\, \cdots \,,\, a_{\,n}\,\right >\,g_{\,i}\;\; \;\forall\; f \,\in\, H_{F}.\]  
\end{theorem}

\begin{proof}
According to the Theorem (3) of \cite{L}, there exists a Bessel sequence \,$\{\,g_{\,i}\,\}_{i \,=\, 1}^{\infty}$\, associated to \,$\left(\,a_{\,2},\, \cdots,\, a_{\,n}\,\right)$\, such that
\[K\,f \,=\, \sum\limits_{i \,=\, 1}^{\infty}\,\left <\,f \,,\,  g_{\,i} \,|\, a_{\,2} \,,\, \cdots \,,\, a_{\,n}\,\right >\,f_{\,i}\;\; \;\forall\; f \,\in\, H_{F}.\]
Now, for each \,$f,\,g \,\in\, H_{F}$, we have
\[\left<\,K\,f \,,\, g \,|\, a_{\,2} \,,\, \cdots \,,\, a_{\,n}\,\right> \,=\, \left<\,\sum\limits_{i \,=\, 1}^{\infty}\,\left<\,f \,,\,  g_{\,i} \,|\, a_{\,2} \,,\, \cdots \,,\, a_{\,n}\,\right>\,f_{\,i} \,,\, g \,|\, a_{\,2} \,,\, \cdots \,,\, a_{\,n}\,\right>\]
\[=\, \sum\limits_{i \,=\, 1}^{\infty}\,\left<\,f \,,\,  g_{\,i} \,|\, a_{\,2} \,,\, \cdots \,,\, a_{\,n}\,\right>\,\left<\,f_{\,i} \,,\,  g \,|\, a_{\,2} \,,\, \cdots \,,\, a_{\,n}\,\right>\hspace{4.1cm}\]
\[=\, \left<\,f \,,\, \sum\limits_{i \,=\, 1}^{\infty}\,\left<\,g \,,\, f_{\,i} \,|\, a_{\,2} \,,\, \cdots \,,\, a_{\,n}\,\right>\,g_{\,i} \,|\, a_{\,2} \,,\, \cdots \,,\, a_{\,n}\,\right>.\hspace{3.6cm}\] 
This shows that \,$K^{\,\ast}\,f \,=\, \sum\limits_{i \,=\, 1}^{\infty}\,\left <\,f \,,\,  f_{\,i} \,|\, a_{\,2} \,,\, \cdots \,,\, a_{\,n}\,\right >\,g_{\,i}\;\; \;\forall\; f \,\in\, H_{F}$.\,This completes the proof.
\end{proof}

\section{Atomic system in Tensor product of $n$-Hilbert spaces}

Let \,$H_{1}$\, and \,$H_{2}$\, be two \,$n$-Hilbert spaces associated with the \,$n$-inner products \,$\left<\,\cdot \,,\, \cdot \,|\, \cdot \,,\, \cdots \,,\, \cdot\,\right>_{1}$\, and \,$\left<\,\cdot \,,\, \cdot \,|\, \cdot \,,\, \cdots \,,\, \cdot\,\right>_{2}$, respectively.\;The tensor product of \,$H_{1}$\, and \,$H_{2}$\, is denoted by \,$H_{1} \,\otimes\, H_{2}$\, and it is defined to be an \,$n$-inner product space associated with the \,$n$-inner product given by
\[\left<\,f_{\,1} \,\otimes\, g_{\,1} \,,\, f_{\,2} \,\otimes\, g_{\,2} \,|\, f_{\,3} \,\otimes\, g_{\,3} \,,\, \,\cdots \,,\, f_{\,n} \,\otimes\, g_{\,n}\,\right>\hspace{1cm}\]
\begin{equation}\label{eqn1}
\hspace{.1cm} \,=\, \left<\,f_{\,1} \,,\, f_{\,2} \,|\, f_{\,3} \,,\, \,\cdots \,,\, f_{\,n}\,\right>_{1}\,\left<\,g_{\,1} \,,\, g_{\,2} \,|\, g_{\,3} \,,\, \,\cdots \,,\, g_{\,n}\,\right>_{2}
\end{equation}
for all \,$f_{\,1},\, f_{\,2},\, f_{\,3},\, \,\cdots,\, f_{\,n} \,\in\, H_{1}$\, and \,$g_{\,1},\, g_{\,2},\, g_{\,3},\, \,\cdots,\, g_{\,n} \,\in\, H_{2}$.\\
The \,$n$-norm on \,$H_{1} \,\otimes\, H_{2}$\, is defined by 
\[\left\|\,f_{\,1} \,\otimes\, g_{\,1},\, f_{\,2} \,\otimes\, g_{\,2},\, \,\cdots,\,\, f_{\,n} \,\otimes\, g_{\,n}\,\right\|\]
\begin{equation}\label{eqn1.1}
 =\,\left\|\,f_{\,1},\, f_{\,2},\, \cdots,\, f_{\,n}\,\right\|_{1}\;\left\|\,g_{\,1},\, g_{\,2},\, \cdots,\, g_{\,n}\,\right\|_{2}
\end{equation}
for all \,$f_{\,1},\, f_{\,2},\, \,\cdots,\, f_{\,n} \,\in\, H_{1}\, \;\text{and}\; \,g_{\,1},\, g_{\,2},\, \,\cdots,\, g_{\,n} \,\in\, H_{2}$, where \,$\left\|\,\cdot \,,\, \cdots \,,\, \cdot \,\right\|_{1}$\, and \,$\left\|\,\cdot \,,\, \cdots \,,\, \cdot \,\right\|_{2}$\, are \,$n$-norm generated by \,$\left<\,\cdot \,,\, \cdot \,|\, \cdot \,,\, \cdots \,,\, \cdot\,\right>_{1}$\, and \,$\left<\,\cdot \,,\, \cdot \,|\, \cdot \,,\, \cdots \,,\, \cdot\,\right>_{2}$, respectively.\;The space \,$H_{1} \,\otimes\, H_{2}$\, is complete with respect to the above \,$n$-inner product.\;Therefore the space \,$H_{1} \,\otimes\, H_{2}$\, is an \,$n$-Hilbert space.     

\begin{note}
Let \,$G \,=\, \left\{\,b_{\,2} \,,\, b_{\,3} \,,\, \cdots \,,\, b_{\,n}\,\right\}$\, be a non-empty finite set, where \,$b_{\,2},\, b_{\,3}$, \,$\cdots,\, b_{\,n}$\, be the fixed elements in \,$H_{2}$.\;Then we define the Hilbert space \,$K_{G}$\, with respect to the inner product is given by
\[\left<\,x \,+\, L_{G}\,,\, y \,+\, L_{G}\,\right>_{G} \,=\, \left<\,x \,,\, y\,\right>_{G} \,=\, \left<\,x \,,\, y \,|\, b_{\,2} \,,\,  \cdots \,,\, b_{\,n} \,\right>_{2}\;\; \;\forall \;\; x,\, y \,\in\, H_{2},\]
where \,$L_{G}$\, denote the linear subspace of \,$H_{2}$\, spanned by the set \,$G$.\;According to the definition (\ref{def0.001}), \,$H_{F} \,\otimes\, K_{G}$\, is the Hilbert space with respect to the inner product: 
\[\left<\,x \,\otimes\, y \,,\, x^{\,\prime} \,\otimes\, y^{\,\prime}\,\right> \,=\, \left<\,x  \,,\, x^{\,\prime}\,\right>_{F}\;\left<\,y  \,,\, y^{\,\prime}\,\right>_{G}\; \;\forall\; x,\, x^{\,\prime} \,\in\, H_{F}\; \;\&\; \;y,\, y^{\,\prime} \,\in\, K_{G}.\]  
\end{note}

\begin{definition}
Let \,$K_{1} \,\in\, \mathcal{B}\,\left(\,H_{F}\,\right)$\, and \,$K_{2} \,\in\, \mathcal{B}\,\left(\,K_{G}\,\right)$.Then the sequence of vectors \,$\left\{\,f_{\,i} \,\otimes\, g_{\,j}\,\right\}_{i,\,j \,=\, 1}^{\,\infty} \,\subseteq\, H_{1} \,\otimes\, H_{2}$\, is said to be a \,$K_{1} \,\otimes\, K_{2}$-frame associated to \,$\left(\,a_{\,2} \,\otimes\, b_{\,2},\, \cdots,\, a_{\,n} \,\otimes\, b_{\,n}\,\right)$\, for \,$H_{1} \,\otimes\, H_{2}$\, if there exist \,$A,\, B \,>\, 0$\, such that
\[A\,\left\|\,\left(\,K_{1} \,\otimes\, K_{2}\,\right)^{\,\ast}\,(\,f \,\otimes\, g\,) \,,\, a_{\,2} \,\otimes\, b_{\,2} \,,\, \cdots \,,\, a_{\,n} \,\otimes\, b_{\,n}\,\right\|^{\,2}\hspace{2cm}\]
\[\,\leq\, \sum\limits_{i,\,j \,=\, 1}^{\,\infty}\,\left|\,\left<\,f \,\otimes\, g,\, f_{\,i} \,\otimes\, g_{\,j} \,|\, a_{\,2} \,\otimes\, b_{\,2},\, \cdots,\, a_{\,n} \,\otimes\, b_{\,n}\,\right>\,\right|^{\,2}\hspace{1.7cm}\]
\begin{equation}\label{em1}
\,\leq\, B\, \left\|\,f \,\otimes\, g,\, a_{\,2} \,\otimes\, b_{\,2},\, \cdots,\, a_{\,n} \,\otimes\, b_{\,n}\,\right\|^{\,2}\; \;\forall\; f \,\otimes\, g \,\in\, H_{F} \,\otimes\, K_{G}.
\end{equation}
If \,$A \,=\, B$, then it is called a tight \,$K_{1} \otimes K_{2}$-frame associated to \,$\left(\,a_{\,2} \otimes b_{\,2},\, \cdots,\, a_{\,n} \otimes b_{\,n}\,\right)$. If \,$K_{1} \,=\, I_{F}$\, and \,$K_{2} \,=\, I_{G}$, then by the Theorem (\ref{th2}), it is a frame associated to \,$\left(\,a_{\,2} \,\otimes\, b_{\,2},\,\cdots,\, a_{\,n} \,\otimes\, b_{\,n}\,\right)$\, for \,$H_{1} \,\otimes\, H_{2}$, where \,$I_{F}$\, and \,$I_{G}$\, are identity operators on \,$H_{F}$\, and \,$K_{G}$, respectively. \\
If only the last inequality of (\ref{em1}) is true then the sequence \,$\left\{\,f_{\,i} \,\otimes\, g_{\,j}\,\right\}_{i,\,j \,=\, 1}^{\,\infty}$\, is called a Bessel sequence associated to \,$\left(\,a_{\,2} \,\otimes\, b_{\,2},\,\cdots,\, a_{\,n} \,\otimes\, b_{\,n}\,\right)$\, in \,$H_{1} \,\otimes\, H_{2}$.\,Thus every  \,$K_{1} \,\otimes\, K_{2}$-frame associated to \,$\left(\,a_{\,2} \otimes b_{\,2},\, \cdots,\, a_{\,n} \otimes b_{\,n}\,\right)$\, is a Bessel sequence associated to \,$\left(\,a_{\,2} \,\otimes\, b_{\,2},\,\cdots,\, a_{\,n} \,\otimes\, b_{\,n}\,\right)$.    
\end{definition}

\begin{theorem}\label{tn3}
Let \,$\left\{\,f_{\,i}\,\right\}_{i \,=\, 1}^{\,\infty}$\, and \,$\left\{\,g_{\,j}\,\right\}_{j \,=\, 1}^{\,\infty}$\, be two sequences in \,$H_{1}$\, and \,$H_{2}$.\;Then \,$\left\{\,f_{\,i}\,\right\}_{i \,=\, 1}^{\,\infty}$\, is a \,$K_{1}$-frame associated to \,$\left(\,a_{\,2},\, \cdots,\, a_{\,n}\,\right)$\, for \,$H_{1}$\, and \,$\left\{\,g_{\,j}\,\right\}_{j \,=\, 1}^{\,\infty}$\, is a \,$K_{2}$-frame associated to \,$\left(\,b_{\,2},\,\cdots,\, b_{\,n}\,\right)$\, for \,$H_{2}$\, if and only if the sequence \,$\left\{\,f_{\,i} \,\otimes\, g_{\,j}\,\right\}_{i,\,j \,=\, 1}^{\,\infty}$\, is a \,$K_{1} \,\otimes\, K_{2}$-frame associated to \,$\left(\,a_{\,2} \,\otimes\, b_{\,2},\, \cdots,\, a_{\,n} \,\otimes\, b_{\,n}\,\right)$\, for \,$H_{1} \,\otimes\, H_{2}$.
\end{theorem}

\begin{proof}
Suppose that the sequence \,$\left\{\,f_{\,i} \,\otimes\, g_{\,j}\,\right\}_{i,\,j \,=\, 1}^{\,\infty}$\, is a \,$K_{1} \,\otimes\, K_{2}$-frame associated to \,$\left(\,a_{\,2} \,\otimes\, b_{\,2},\, \cdots,\, a_{\,n} \,\otimes\, b_{\,n}\,\right)$\, for \,$H_{1} \,\otimes\, H_{2}$.\;Then for each \,$f \,\otimes\, g \,\in\, H_{F} \,\otimes\, K_{G} \,-\, \{\,\theta \,\otimes\, \theta\,\}$, there exist constants \,$A,\,B \,>\, 0$\, such that
\[A\,\left\|\,\left(\,K_{1} \,\otimes\, K_{2}\,\right)^{\,\ast}\,(\,f \,\otimes\, g\,) \,,\, a_{\,2} \,\otimes\, b_{\,2} \,,\, \cdots \,,\, a_{\,n} \,\otimes\, b_{\,n}\,\right\|^{\,2}\]
\[\,\leq\, \sum\limits_{i,\,j \,=\, 1}^{\,\infty}\,\left|\,\left<\,f \,\otimes\, g,\, f_{\,i} \,\otimes\, g_{\,j} \,|\, a_{\,2} \,\otimes\, b_{\,2},\, \cdots,\, a_{\,n} \,\otimes\, b_{\,n}\,\right>\,\right|^{\,2}\]
\[ \,\leq\, B\, \left\|\,f \,\otimes\, g,\, a_{\,2} \,\otimes\, b_{\,2},\, \cdots,\, a_{\,n} \,\otimes\, b_{\,n}\,\right\|^{\,2}\hspace{2.4cm}\]
\[\Rightarrow\,A\left\|\,K^{\,\ast}_{1}\,f \,\otimes\, K^{\,\ast}_{2}\,g,\, a_{\,2} \,\otimes\, b_{\,2},\, \cdots,\, a_{\,n} \,\otimes\, b_{\,n}\,\right\|^{\,2}\hspace{6.1cm}\]
\[\,\leq\, \sum\limits_{i,\,j \,=\, 1}^{\,\infty}\,\left|\,\left<\,f \,\otimes\, g,\, f_{\,i} \,\otimes\, g_{\,j} \,|\, a_{\,2} \,\otimes\, b_{\,2},\, \cdots,\, a_{\,n} \,\otimes\, b_{\,n}\,\right>\,\right|^{\,2}\]
\[ \,\leq\, B\, \left\|\,f \,\otimes\, g,\, a_{\,2} \,\otimes\, b_{\,2},\, \cdots,\, a_{\,n} \,\otimes\, b_{\,n}\,\right\|^{\,2}\hspace{2.2cm}\]
\[\Rightarrow A\left\|\,K_{1}^{\,\ast}f,\, a_{2},\, \cdots,\, a_{n}\,\right\|_{1}^{\,2}\,\left\|\,K_{2}^{\,\ast}g,\, b_{2},\, \cdots,\, b_{n}\,\right\|_{2}^{\,2} \leq \left(\sum\limits_{i \,=\, 1}^{\,\infty}\left|\,\left<\,f,\, f_{\,i}\,|\,a_{2},\, \cdots,\, a_{n}\,\right>_{1}\,\right|^{\,2}\right)\times\]
\[\left(\,\sum\limits_{j \,=\, 1}^{\,\infty}\,\left|\,\left<\,g \,,\, g_{\,j} \,|\, b_{\,2} \,,\, \cdots \,,\, b_{\,n}\,\right>_{2}\,\right|^{\,2}\,\right) \,\leq\, B\,\left\|\,f,\, a_{\,2},\, \cdots,\, a_{\,n}\,\right\|_{1}^{\,2}\,\left\|\,g,\, b_{\,2},\, \cdots,\, b_{\,n}\,\right\|_{2}^{\,2}.\]
Since \,$f \,\otimes\, g \,\in\, H_{F} \,\otimes\, K_{G}$\, is non-zero, \,$f \,\in\, H_{F}$\, and \,$g \,\in\, K_{G}$\, are non-zero elements and therefore \,$\sum\limits_{\,j \,=\, 1}^{\,\infty}\,\left|\,\left<\,g \,,\, g_{\,j} \,|\, b_{\,2} \,,\, \cdots \,,\, b_{\,n}\,\right>_{2}\,\right|^{\,2}$\, and \,$\sum\limits_{\,i \,=\, 1}^{\,\infty}\,\left|\,\left<\,f \,,\, f_{\,i} \,|\, a_{\,2} \,,\, \cdots \,,\, a_{\,n}\,\right>_{1}\,\right|^{\,2}$\, are non-zero.\,This implies that
\[\dfrac{A\left\|\,K_{2}^{\,\ast}g,\, b_{\,2},\, \cdots,\, b_{\,n}\,\right\|_{2}^{\,2}}{\sum\limits_{j \,=\, 1}^{\,\infty}\,\left|\,\left<\,g,\, g_{\,j}\,|\,b_{\,2},\, \cdots,\, b_{\,n}\,\right>_{2}\,\right|^{\,2}}\,\left\|\,K_{1}^{\,\ast}\,f,\, a_{\,2},\, \cdots,\, a_{\,n}\,\right\|_{1}^{\,2} \,\leq\, \sum\limits_{i \,=\, 1}^{\,\infty}\,\left|\,\left<\,f,\, f_{\,i}\,|\,a_{\,2},\, \cdots,\, a_{\,n}\,\right>_{1}\,\right|^{\,2}\]
\[\hspace{2cm}\leq\,  \dfrac{B\,\left\|\,K_{2}^{\,\ast}\,g,\, b_{\,2},\, \cdots,\, b_{\,n}\,\right\|_{2}^{\,2}}{\sum\limits_{j \,=\, 1}^{\,\infty}\,\left|\,\left<\,g \,,\, g_{\,j} \,|\, b_{\,2} \,,\, \cdots \,,\, b_{\,n}\,\right>_{2}\,\right|^{\,2}}\,\left\|\,f,\, a_{\,2},\, \cdots,\, a_{\,n}\,\right\|_{1}^{\,2}\]
\[\Rightarrow A_{1}\left\|\,K_{1}^{\,\ast}\,f,\, a_{\,2},\, \cdots,\, a_{\,n}\,\right\|_{1}^{\,2} \leq \sum\limits_{i \,=\, 1}^{\,\infty}\left|\,\left<\,f,\, f_{\,i}\,|\,a_{\,2},\, \cdots,\, a_{\,n}\,\right>_{1}\,\right|^{\,2} \leq B_{1}\,\left\|\,f,\, a_{\,2},\, \cdots,\, a_{\,n}\,\right\|_{1}^{\,2},\]
where \,$A_{1} \,=\, \dfrac{A\,\left\|\,K_{2}^{\,\ast}\,g,\, b_{\,2},\, \cdots,\, b_{\,n}\,\right\|_{2}^{\,2}}{\sum\limits_{j \,=\, 1}^{\,\infty}\,\left|\,\left<\,g,\, g_{\,j}\,|\,b_{\,2},\, \cdots,\, b_{\,n}\,\right>_{2}\,\right|^{\,2}}$\, and \,$B_{1} \,=\, \dfrac{B\,\left\|\,K_{2}^{\,\ast}\,g,\, b_{\,2},\, \cdots,\, b_{\,n}\,\right\|_{2}^{\,2}}{\sum\limits_{j \,=\, 1}^{\,\infty}\,\left|\,\left<\,g,\, g_{\,j}\,|\,b_{\,2},\, \cdots,\, b_{\,n}\,\right>_{2}\,\right|^{\,2}}$. This shows that \,$\{\,f_{\,i}\,\}_{i \,=\,1}^{\infty}$\, is a \,$K_{1}$-frame associated to \,$\left(\,a_{\,2},\, \cdots,\, a_{\,n}\,\right)$\, for \,$H_{1}$.\;Similarly, it can be shown that \,$\{\,g_{\,j}\,\}_{j \,=\,1}^{\infty}$\, is a \,$K_{2}$-frame associated to \,$\left(\,b_{\,2},\, \cdots,\, b_{\,n}\,\right)$\, for \,$H_{2}$.\\

Conversely, Suppose that \,$\left\{\,f_{\,i}\,\right\}_{i \,=\, 1}^{\,\infty}$\, is a \,$K_{1}$-frame associated to \,$\left(\,a_{\,2},\, \cdots,\, a_{\,n}\,\right)$\, for \,$H_{1}$\, with bounds \,$A,\, B$\, and \,$\left\{\,g_{\,j}\,\right\}_{j \,=\, 1}^{\,\infty}$\, is a \,$K_{2}$-frame associated to \,$\left(\,b_{\,2},\, \cdots,\, b_{\,n}\,\right)$\, for \,$H_{2}$\, with bounds \,$C,\, D$.\;Then, for all \,$f \,\in\, H_{F}$\, and \,$g \,\in\, K_{G}$, we have
\[A\left\|\,K_{1}^{\,\ast}\,f,\, a_{\,2},\, \cdots,\, a_{\,n}\,\right\|_{1}^{\,2} \leq \sum\limits_{i \,=\, 1}^{\,\infty}\,\left|\,\left<\,f,\, f_{\,i}\,|\,a_{\,2},\, \cdots,\, a_{\,n}\,\right>_{1}\,\right|^{\,2} \,\leq\, B\,\left\|\,f,\, a_{\,2},\, \cdots,\, a_{\,n}\,\right\|_{1}^{\,2},\;\&\]
\[C\,\left\|\,K_{2}^{\,\ast}\,g,\, b_{\,2},\, \cdots,\, b_{\,n}\,\right\|_{2}^{\,2} \,\leq\, \sum\limits_{j \,=\, 1}^{\,\infty}\,\left|\,\left<\,g,\, g_{\,j} \,|\, b_{\,2},\, \cdots,\, b_{\,n}\,\right>_{2}\,\right|^{\,2} \,\leq\, D\,\left\|\,g,\, b_{\,2},\, \cdots,\, b_{\,n}\,\right\|_{2}^{\,2}.\]
Multiplying the above two inequalities and using (\ref{eqn1}) and (\ref{eqn1.1}), we get
\[A\,C\,\left\|\,\left(\,K_{1} \,\otimes\, K_{2}\,\right)^{\,\ast}\,(\,f \,\otimes\, g\,) \,,\, a_{\,2} \,\otimes\, b_{\,2} \,,\, \cdots \,,\, a_{\,n} \,\otimes\, b_{\,n}\,\right\|^{\,2}\hspace{1cm}\]
\[ \,\leq\, \sum\limits_{i,\,j \,=\, 1}^{\,\infty}\,\left|\,\left<\,f \,\otimes\, g,\, f_{\,i} \,\otimes\, g_{\,j} \,|\, a_{\,2} \,\otimes\, b_{\,2}, \,\cdots,\, a_{\,n} \,\otimes\, b_{\,n}\,\right>\,\right|^{\,2}\hspace{1cm}\]
\[\hspace{1cm}\leq\, B\,D\,\left\|\,f \,\otimes\, g,\, a_{\,2} \,\otimes\, b_{\,2}, \,\cdots,\, a_{\,n} \,\otimes\, b_{\,n}\,\right\|^{\,2}\; \;\forall\; f \,\otimes\, g \,\in\, H_{F} \,\otimes\, K_{G}.\]  
Hence, \,$\left\{\,f_{\,i} \,\otimes\, g_{\,j}\,\right\}_{i,\,j \,=\, 1}^{\,\infty}$\, is a \,$K_{1} \,\otimes\, K_{2}$-frame associated to \,$\left(\,a_{\,2} \,\otimes\, b_{\,2},\, \cdots,\, a_{\,n} \,\otimes\, b_{\,n}\,\right)$\, for \,$H_{1} \,\otimes\, H_{2}$.\;This completes the proof.
\end{proof}

\begin{theorem}\label{th3.01}
Let \,$\{\,f_{\,i}\,\}_{i \,=\,1}^{\infty}$\, and \,$\{\,g_{\,j}\,\}_{j \,=\,1}^{\infty}$\, be the sequences of vectors in \,$n$-Hilbert spaces \,$H_{1}$\, and \,$H_{2}$.\;Then the sequence \,$\left\{\,f_{\,i} \,\otimes\, g_{\,j}\,\right\}^{\,\infty}_{i,\,j \,=\, 1} \,\subseteq\, H_{1} \,\otimes\, H_{2}$\, is a Bessel sequence associated to \,$\left(\,a_{\,2} \,\otimes\, b_{\,2},\, \,\cdots,\, a_{\,n} \,\otimes\, b_{\,n}\,\right)$\, in \,$H_{1} \,\otimes\, H_{2}$\, if and only if \,$\{\,f_{\,i}\,\}_{i \,=\,1}^{\infty}$\, is a Bessel sequence associated to \,$\left(\,a_{\,2},\, \cdots,\, a_{\,n}\,\right)$\, in \,$H_{1}$\, and \,$\{\,g_{\,j}\,\}_{j \,=\,1}^{\infty}$\, is a Bessel sequence associated to \,$\left(\,b_{\,2},\, \cdots,\, b_{\,n}\,\right)$\, in \,$H_{2}$.   
\end{theorem}

\begin{proof}
Since every  \,$K_{1} \,\otimes\, K_{2}$-frame associated to \,$\left(\,a_{\,2} \otimes b_{\,2},\, \cdots,\, a_{\,n} \otimes b_{\,n}\,\right)$\, is a Bessel sequence associated to \,$\left(\,a_{\,2} \,\otimes\, b_{\,2},\,\cdots,\, a_{\,n} \,\otimes\, b_{\,n}\,\right)$, proof of this theorem directly follows from the Theorem (\ref{tn3}).
\end{proof}

\begin{theorem}
Let \,$\left\{\,f_{\,i}\,\right\}_{i \,=\, 1}^{\,\infty}$\, be a \,$K_{1}$-frame associated to \,$\left(\,a_{\,2},\, \cdots,\, a_{\,n}\,\right)$\, for \,$H_{1}$\, with bounds \,$A,\,B$\, and \,$\left\{\,g_{\,j}\,\right\}_{j \,=\, 1}^{\,\infty}$\, be a \,$K_{2}$-frame associated to \,$\left(\,b_{\,2},\, \cdots,\, b_{\,n}\,\right)$\, for \,$H_{2}$\, with bounds \,$C,\,D$, respectively.

\begin{itemize}
\item[(I)]If \,$T_{1} \,\otimes\, T_{2} \,\in\, \mathcal{B}\,\left(\,H_{F} \,\otimes\, K_{G}\,\right)$\, is an isometry such that \,$\left(\,K_{1} \otimes K_{2}\,\right)^{\,\ast}\left(\,T_{1} \otimes T_{2}\,\right) = \left(\,T_{1} \otimes T_{2}\,\right)\,\left(\,K_{1} \otimes K_{2}\,\right)^{\,\ast}$, then the sequence \,$\Delta \,=\, \left\{\,\left(\,T_{1} \otimes T_{2}\,\right)^{\,\ast}\,\left(\,f_{\,i} \,\otimes\, g_{\,j}\,\right)\,\right\}_{i,\,j \,=\, 1}^{\,\infty}$\, is a \,$K_{1} \,\otimes\, K_{2}$-frame associated to \,$\left(\,a_{\,2} \,\otimes\, b_{\,2},\, \cdots,\, a_{\,n} \,\otimes\, b_{\,n}\,\right)$\, for \,$H_{1} \,\otimes\, H_{2}$.
\item[(II)]The sequence \,$\Gamma \,=\, \left\{\,\left(\,L_{1} \,\otimes\, L_{2}\,\right)\,\left(\,f_{\,i} \,\otimes\, g_{\,j}\,\right)\,\right\}_{i,\,j \,=\, 1}^{\,\infty}$\, is a \,$\left(\,L_{1} \,\otimes\, L_{2}\,\right)\,\left(\,K_{1} \,\otimes\, K_{2}\,\right)$-frame associated to \,$\left(\,a_{\,2} \,\otimes\, b_{\,2},\, \cdots,\, a_{\,n} \,\otimes\, b_{\,n}\,\right)$\, for \,$H_{1} \,\otimes\, H_{2}$, for some operator \,$L_{1} \,\otimes\, L_{2} \,\in\, \mathcal{B}\,\left(\,H_{F} \,\otimes\, K_{G}\,\right)$.  
\end{itemize} 
\end{theorem}

\begin{proof}$(I)$\,
For each \,$f \,\otimes\, g \,\in\, H_{F} \,\otimes\, K_{G}$, we have
\[\sum\limits_{i,\,j \,=\, 1}^{\,\infty}\,\left|\,\left<\,f \,\otimes\, g \,,\, \left(\,T_{1} \,\otimes\, T_{2}\,\right)^{\,\ast}\,\left(\,f_{\,i} \,\otimes\, g_{\,j}\,\right) \,|\, a_{\,2} \,\otimes\, b_{\,2},\, \cdots,\, a_{\,n} \,\otimes\, b_{\,n}\,\right>\,\right|^{\,2}\]
\[ \,=\,\sum\limits_{i,\,j \,=\, 1}^{\,\infty}\,\left|\,\left<\,f \,\otimes\, g \,,\, T^{\,\ast}_{1}\,f_{\,i} \,\otimes\, T^{\,\ast}_{2}\,g_{\,j} \,|\, a_{\,2} \,\otimes\, b_{\,2},\, \cdots,\, a_{\,n} \,\otimes\, b_{\,n}\,\right>\,\right|^{\,2}\hspace{3.2cm}\]
\[=\, \left(\,\sum\limits_{i \,=\, 1}^{\,\infty}\,\left|\,\left<\,f \,,\, T^{\,\ast}_{1}\,f_{\,i} \,|\, a_{\,2},\, \cdots,\, a_{\,n}\,\right>_{1}\,\right|^{\,2}\,\right)\,\left(\,\sum\limits_{j \,=\, 1}^{\,\infty}\left|\,\left<\,g \,,\, T^{\,\ast}_{2}\,g_{\,j} \,|\, b_{\,2},\, \cdots,\, b_{\,n}\,\right>_{2}\,\right|^{\,2}\,\right)\hspace{.3cm}\]
\begin{equation}\label{eq1}
=\, \left(\,\sum\limits_{i \,=\, 1}^{\,\infty}\,\left|\,\left<\,T_{1}\,f \,,\, f_{\,i} \,|\, a_{\,2},\, \cdots,\, a_{\,n}\,\right>_{1}\,\right|^{\,2}\,\right)\,\left(\,\sum\limits_{j \,=\, 1}^{\,\infty}\left|\,\left<\,T_{2}\,g \,,\, g_{\,j} \,|\, b_{\,2},\, \cdots,\, b_{\,n}\,\right>_{2}\,\right|^{\,2}\,\right)
\end{equation}
\[\leq\, B\,\left\|\,T_{1}\,f \,,\, a_{\,2},\, \cdots,\, a_{\,n}\,\right\|_{1}^{\,2}\,D\,\left\|\,T_{2}\,g \,,\, a_{\,2},\, \cdots,\, a_{\,n}\,\right\|_{2}^{\,2}\hspace{4.7cm}\]
\[[\;\text{since \,$\left\{\,f_{\,i}\,\right\}_{i \,=\, 1}^{\,\infty}$\, is a \,$K_{1}$-frame associated to \,$\left(\,a_{\,2},\, \cdots,\, a_{\,n}\,\right)$\, and}\]
\[\text{ \,$\left\{\,g_{\,j}\,\right\}_{j \,=\, 1}^{\,\infty}$\, is a \,$K_{2}$-frame associated to \,$\left(\,b_{\,2},\, \cdots,\, b_{\,n}\,\right)$}\;]\]
\[\leq\, B\,D\,\left\|\,T_{1}\,\right\|^{2}\,\left\|\,T_{2}\,\right\|^{2}\,\left\|\,f \,,\, a_{\,2},\, \cdots,\, a_{\,n}\,\right\|^{\,2}_{1}\,\left\|\,g \,,\, b_{\,2},\, \cdots,\, b_{\,n}\,\right\|^{\,2}_{2}\hspace{3.3cm}\]
\[ \,=\, B\,D\,\left\|\,T_{1} \,\otimes\, T_{2}\,\right\|^{\,2}\,\left\|\,f \,\otimes\, g \,,\, a_{\,2} \,\otimes\, b_{\,2},\, \cdots,\, a_{\,n} \,\otimes\, b_{\,n}\,\right\|^{\,2}.\hspace{3.7cm}\]
On the other hand, since \,$\left\{\,f_{\,i}\,\right\}_{i \,=\, 1}^{\,\infty}$\, is a \,$K_{1}$-frame associated to \,$\left(\,a_{\,2},\, \cdots,\, a_{\,n}\,\right)$\, for \,$H_{1}$\, and \,$\left\{\,g_{\,j}\,\right\}_{j \,=\, 1}^{\,\infty}$\, is a \,$K_{2}$-frame associated to \,$\left(\,b_{\,2},\, \cdots,\, b_{\,n}\,\right)$\, for \,$H_{2}$, from (\ref{eq1}),
\[\sum\limits_{i,\,j \,=\, 1}^{\,\infty}\,\left|\,\left<\,f \,\otimes\, g \,,\, \left(\,T_{1} \,\otimes\, T_{2}\,\right)^{\,\ast}\,\left(\,f_{\,i} \,\otimes\, g_{\,j}\,\right) \,|\, a_{\,2} \,\otimes\, b_{\,2},\, \cdots,\, a_{\,n} \,\otimes\, b_{\,n}\,\right>\,\right|^{\,2}\hspace{2.5cm}\]
\[ \,\geq\, A\,\left\|\,K^{\,\ast}_{1}\,T_{1}\,f \,,\, a_{\,2},\, \cdots,\, a_{\,n}\,\right\|_{1}^{\,2}\,C\,\left\|\,K^{\,\ast}_{2}\,T_{2}\,g \,,\, b_{\,2},\, \cdots,\, b_{\,n}\,\right\|_{2}^{\,2}\hspace{4.5cm}\]
\[ \,=\, A\, C\,\left<\,K^{\,\ast}_{1}\,T_{\,1}\,f \,,\, K^{\,\ast}_{1}\,T_{1}\,f \;|\; a_{\,2},\, \cdots,\, a_{\,n}\,\right>_{\,1}\;\left<\,K^{\,\ast}_{2}\,T_{2}\,g \,,\, K^{\,\ast}_{2}\,T_{2}\,g \;|\; b_{\,2},\, \cdots,\, b_{\,n}\,\right>_{\,2}\hspace{2.1cm}\]
\[=\, A\,C\,\left<\,K^{\,\ast}_{\,1}\,T_{\,1}\,f \,\otimes\, K^{\,\ast}_{\,2}\,T_{\,2}\,g \;,\; K^{\,\ast}_{\,1}\,T_{\,1}\,f \,\otimes\, K^{\,\ast}_{\,2}\,T_{\,2}\,g \;|\; a_{\,2} \,\otimes\, b_{\,2},\, \cdots,\, a_{\,n} \,\otimes\, b_{\,n}\,\right>\hspace{4.5cm}\]
\[= A\,C\left<\left(K_{1}\otimes K_{2}\right)^{\,\ast}\left(T_{1} \otimes T_{2}\right)(\,f \otimes g\,),\, \left(K_{1} \otimes K_{2}\right)^{\,\ast}\left(T_{1} \otimes T_{2}\right)(\,f \otimes g\,)\,|\,a_{2} \otimes b_{2},\, \cdots,\, a_{n} \otimes b_{n}\right>\]
\[= A\,C\left<\left(T_{1} \otimes T_{2}\right)\left(K_{1} \otimes K_{2}\right)^{\,\ast}(\,f \otimes g\,),\, \left(T_{1} \otimes T_{2}\right)\left(K_{1}\otimes K_{2}\right)^{\,\ast}(\,f \otimes g\,)\,|\,a_{2} \otimes b_{2},\, \cdots,\, a_{n} \otimes b_{n}\right>\]
\[\hspace{2cm}\left[\,\text{since}\; \left(\,K_{1} \,\otimes\, K_{2}\,\right)^{\,\ast}\,\left(\,T_{1} \,\otimes\, T_{2}\,\right) \,=\, \left(\,T_{1} \,\otimes\, T_{2}\,\right)\,\left(\,K_{1} \,\otimes\, K_{2}\,\right)^{\,\ast}\,\right]\]
\[=\, A\,C\,\left<\,\left(\,K_{1} \,\otimes\, K_{2}\,\right)^{\,\ast}\,(\,f \,\otimes\, g\,) \;,\; \left(\,K_{1} \,\otimes\, K_{2}\,\right)^{\,\ast}\,(\,f \,\otimes\, g\,) \;|\; a_{\,2} \,\otimes\, b_{\,2},\, \cdots,\, a_{\,n} \,\otimes\, b_{\,n}\,\right>\]
\[\hspace{2cm}[\;\text{since}\; \left(\,T_{1} \,\otimes\, T_{2}\,\right)\; \;\text{is an isometry}\;]\]  
\[=\, A\,C\,\left<\,K^{\,\ast}_{1}\,f \,\otimes\, K^{\,\ast}_{2}\,g \;,\; K^{\,\ast}_{1}\,f \,\otimes\, K^{\,\ast}_{2}\,g \;|\; a_{\,2} \,\otimes\, b_{\,2},\, \cdots,\, a_{\,n} \,\otimes\, b_{\,n}\,\right>\hspace{5cm}\]
\[ \,=\, A\,C\,\left<\,K^{\,\ast}_{1}\,f \,,\, K^{\,\ast}_{1}\,f \,|\, a_{\,2},\, \cdots,\, a_{\,n}\,\right>_{\,1}\,\left<\,K^{\,\ast}_{2}\,g \,,\, K^{\,\ast}_{2}\,g \,|\, b_{\,2},\, \cdots,\, b_{\,n}\,\right>_{\,2}\hspace{5cm}\]
\[=\, A\,C\,\left\|\,K^{\,\ast}_{1}\,f \,,\, a_{\,2},\, \cdots,\, a_{\,n}\,\right\|_{\,1}^{\,2}\,\left\|\,K^{\,\ast}_{2}\,g \,,\, b_{\,2},\, \cdots,\, b_{\,n}\,\right\|_{\,2}^{\,2}\hspace{5cm}\]
\[ \,=\, A\,C\,\left\|\,K^{\,\ast}_{1}\,f \,\otimes\, K^{\,\ast}_{2}\,g \,,\, a_{\,2} \,\otimes\, b_{\,2},\, \cdots,\, a_{\,n} \,\otimes\, b_{\,n}\,\right\|^{\,2}\hspace{5.5cm}\]
\[=\, A\,C\,\left\|\,\left(\,K_{1} \,\otimes\, K_{2}\,\right)^{\,\ast}\,(\,f \,\otimes\, g\,) \,,\, a_{\,2} \,\otimes\, b_{\,2},\, \cdots,\, a_{\,n} \,\otimes\, b_{\,n}\,\right\|^{\,2}.\hspace{5.5cm}\]
Hence, \,$\Delta$\, is a \,$K_{1} \,\otimes\, K_{2}$-frame associated to \,$\left(\,a_{\,2} \,\otimes\, b_{\,2},\, \cdots,\, a_{\,n} \,\otimes\, b_{\,n}\,\right)$\, for \,$H_{1} \,\otimes\, H_{2}$.\\

$(II)$\, According to the proof of $(I)$, it is easy to verify that
\[\sum\limits_{i,\,j \,=\, 1}^{\,\infty}\,\left|\,\left<\,f \,\otimes\, g \,,\, \left(\,L_{1} \,\otimes\, L_{2}\,\right)\,\left(\,f_{\,i} \,\otimes\, g_{\,j}\,\right) \,|\, a_{\,2} \,\otimes\, b_{\,2},\, \cdots,\, a_{\,n} \,\otimes\, b_{\,n}\,\right>\,\right|^{\,2}\hspace{1.5cm}\]
\[ \,\leq\, B\,D\,\left\|\,L_{1} \otimes L_{2}\,\right\|^{\,2}\left\|\,f \otimes g,\, a_{\,2} \otimes b_{\,2},\, \cdots,\, a_{\,n} \otimes b_{\,n}\,\right\|^{\,2}\, \,\forall\, f \otimes g \,\in\, H_{F} \otimes K_{G}.\]
On the other hand,
\[\sum\limits_{i,\,j \,=\, 1}^{\,\infty}\,\left|\,\left<\,f \,\otimes\, g \,,\, \left(\,L_{1} \,\otimes\, L_{2}\,\right)\,\left(\,f_{\,i} \,\otimes\, g_{\,j}\,\right) \,|\, a_{\,2} \,\otimes\, b_{\,2},\, \cdots,\, a_{\,n} \,\otimes\, b_{\,n}\,\right>\,\right|^{\,2}\hspace{1.2cm}\]
\[ \,\geq\, A\,\left\|\,K^{\,\ast}_{1}\,L^{\,\ast}_{1}\,f \,,\, a_{\,2},\, \cdots,\, a_{\,n}\,\right\|_{1}^{\,2}\,C\,\left\|\,K^{\,\ast}_{2}\,L^{\,\ast}_{2}\,g \,,\, b_{\,2},\, \cdots,\, b_{\,n}\,\right\|_{2}^{\,2}\hspace{2.4cm}\]
\[ \,=\, A\,C\,\left\|\,K^{\,\ast}_{1}\,L^{\,\ast}_{1}\,f \,\otimes\, K^{\,\ast}_{2}\,L^{\,\ast}_{2}\,g \,,\, a_{\,2} \,\otimes\, b_{\,2},\, \cdots,\, a_{\,n} \,\otimes\, b_{\,n}\,\right\|^{\,2}\hspace{2.9cm}\]
\[=\, A\,C\,\left\|\,\left(\,K^{\,\ast}_{1}\,L^{\,\ast}_{1} \,\otimes\, K^{\,\ast}_{2}\,L^{\,\ast}_{2}\,\right)\,(\,f \,\otimes\, g\,) \,,\, a_{\,2} \,\otimes\, b_{\,2},\, \cdots,\, a_{\,n} \,\otimes\, b_{\,n}\,\right\|^{\,2}\hspace{1.4cm}\]
\[ \,=\, A\,C\,\left\|\,\left[\,\left(\,L_{1} \,\otimes\, L_{2}\,\right)\,\left(\,K_{1} \,\otimes\, K_{2}\,\right)\,\right]^{\,\ast}\,(\,f \,\otimes\, g\,) \,,\, a_{\,2} \,\otimes\, b_{\,2},\, \cdots,\, a_{\,n} \,\otimes\, b_{\,n}\,\right\|^{\,2}.\]
Hence, \,$\Gamma$\, is a \,$\left(\,L_{1} \,\otimes\, L_{2}\,\right)\,\left(\,K_{1} \,\otimes\, K_{2}\,\right)$-frame associated to \,$\left(\,a_{\,2} \,\otimes\, b_{\,2},\, \cdots,\, a_{\,n} \,\otimes\, b_{\,n}\,\right)$\, for \,$H_{1} \,\otimes\, H_{2}$.
\end{proof}

\begin{definition}\label{def1}
Let \,$K_{1}$\, and \,$K_{2}$\, be bounded linear operators on the Hilbert spaces \,$H_{F}$\, and \,$K_{G}$.\;Then the sequence of vectors \,$\left\{\,f_{\,i} \,\otimes\, g_{\,j}\,\right\}_{i,\,j \,=\, 1}^{\,\infty} \,\subseteq\, H_{1} \,\otimes\, H_{2}$\, is said to be an atomic system associated to \,$\left(\,a_{\,2} \,\otimes\, b_{\,2},\, \cdots,\, a_{\,n} \,\otimes\, b_{\,n}\,\right)$\, for \,$K_{1} \,\otimes\, K_{2} \,\in\, \mathcal{B}\,\left(\,H_{F} \,\otimes\, K_{G}\,\right)$\, in \,$H_{1} \,\otimes\, H_{2}$\, if 
\begin{itemize}
\item[(I)]\,$\left\{\,f_{\,i} \,\otimes\, g_{\,j}\,\right\}_{i,\,j \,=\, 1}^{\,\infty}$\, is a Bessel sequence associated to \,$\left(\,a_{\,2} \,\otimes\, b_{\,2},\, \cdots,\, a_{\,n} \,\otimes\, b_{\,n}\,\right)$\, in \,$H_{1} \,\otimes\, H_{2}$.
\item[(II)] For any \,$f \,\otimes\, g \,\in\, H_{F} \,\otimes\, K_{G}$, there exists \,$c \,\otimes\, d \,=\, \{\,c_{\,i}\,d_{\,j}\,\}_{i,\,j \,=\, 1}^{\,\infty} \,\in\, l^{\,2}\,(\,\mathbb{N} \,\times\, \mathbb{N}\,)$\, such that \,$\left(\,K_{1} \,\otimes\, K_{2}\,\right)\,(\,f \,\otimes\, g\,) \,=\, \sum\limits_{\,i,\, j \,=\, 1}^{\,\infty}\,c_{\,i}\,d_{\,j}\,\left(\,f_{\,i} \,\otimes\, g_{\,j}\,\right)$, and for some \,$C \,>\, 0$, \,$\left\|\,c \,\otimes\, d\,\right\|_{l^{2}} \,\leq\, C\, \left\|\,f \otimes g,\, a_{\,2} \otimes b_{\,2},\, \cdots,\, a_{\,n} \otimes b_{\,n}\,\right\|$, where \,$c \,=\, \{\,c_{\,i}\,\}_{i \,=\,1}^{\,\infty}$, \,$d \,=\, \{\,d_{\,j}\,\}_{j \,=\,1}^{\,\infty}$\, are in \,$l^{\,2}\,(\,\mathbb{N}\,)$.  
\end{itemize} 
\end{definition}

\begin{theorem}\label{th3.1}
Let \,$\left\{\,f_{\,i}\,\right\}_{i \,=\, 1}^{\,\infty}$\, be an atomic system associated to \,$\left(\,a_{\,2},\, \cdots,\, a_{\,n}\,\right)$\, for \,$K_{1}$\, in \,$H_{1}$\, and \,$\left\{\,g_{\,j}\,\right\}_{j \,=\, 1}^{\,\infty}$\, be an atomic system associated to \,$\left(\,b_{\,2},\, \cdots,\, b_{\,n}\,\right)$\, for \,$K_{2}$\, in \,$H_{2}$.\;Then the sequence \,$\left\{\,f_{\,i} \,\otimes\, g_{\,j}\,\right\}_{i,\,j \,=\, 1}^{\,\infty}$\, is an atomic system  associated to \,$\left(\,a_{\,2} \,\otimes\, b_{\,2},\, \cdots,\, a_{\,n} \,\otimes\, b_{\,n}\,\right)$\, for \,$K_{1} \,\otimes\, K_{2}$\, in \,$H_{1} \,\otimes\, H_{2}$.
\end{theorem}

\begin{proof}
Since \,$\left\{\,f_{\,i}\,\right\}_{i \,=\, 1}^{\,\infty}$\, is an atomic system associated to \,$\left(\,a_{\,2},\, \cdots,\, a_{\,n}\,\right)$\, for \,$K_{1}$\, in \,$H_{1}$\, and \,$\left\{\,g_{\,j}\,\right\}_{j \,=\, 1}^{\,\infty}$\, is an atomic system associated to \,$\left(\,b_{\,2},\, \cdots,\, b_{\,n}\,\right)$\, for \,$K_{2}$\, in \,$H_{2}$, by the definition (\ref{defn1}),  \,$\left\{\,f_{\,i}\,\right\}_{i \,=\, 1}^{\,\infty}$\, is a Bessel sequence associated to \,$\left(\,a_{\,2},\, \cdots,\, a_{\,n}\,\right)$\, in \,$H_{1}$\, and \,$\left\{\,g_{\,j}\,\right\}_{j \,=\, 1}^{\,\infty}$\, is a Bessel sequence associated to \,$\left(\,b_{\,2},\, \cdots,\, b_{\,n}\,\right)$\, in \,$H_{2}$, respectively.\;Then by the Theorem (\ref{th3.01}), \,$\left\{\,f_{\,i} \,\otimes\, g_{\,j}\,\right\}_{i,\,j \,=\, 1}^{\,\infty}$\, is a Bessel sequence associated to \,$\left(\,a_{\,2} \,\otimes\, b_{\,2},\, \cdots,\, a_{\,n} \,\otimes\, b_{\,n}\,\right)$\, in \,$H_{1} \,\otimes\, H_{2}$.\;Also, for any \,$f \,\in\, H_{F}$\, and \,$g \,\in\, K_{G}$, 
\[K_{1}\,f \,= \sum\limits_{i \,=\, 1}^{\,\infty}\,c_{\,i}\,f_{\,i}\; \;\text{with}\; \left\|\,\left\{\,c_{\,i}\,\right\}_{i \,=\, 1}^{\,\infty}\,\right\|_{l^{\,2}} \,\leq\, C_{1} \left\|\,f,\, a_{\,2},\, \cdots,\, a_{\,n}\,\right\|_{\,1},\,\text{for some \,$C_{1} \,>\, 0$}\] 
\[K_{2}\,g \,= \sum\limits_{j \,=\, 1}^{\,\infty}\,d_{\,j}\,g_{\,j}\; \,\text{with}\; \left\|\,\left\{\,d_{\,j}\,\right\}_{j \,=\, 1}^{\,\infty}\,\right\|_{l^{\,2}} \leq\, C_{2} \left\|\,g,\, b_{\,2},\, \cdots,\, b_{\,n}\,\right\|_{\,2},\,\text{for some \,$C_{2} \,>\, 0$}.\]Therefore, for each \,$f \,\otimes\, g \,\in\, H_{F} \,\otimes\, K_{G}$, we have 
\[\left(\,K_{1} \otimes K_{2}\,\right)(\,f \,\otimes\, g\,) \,=\, K_{1}f \,\otimes\, K_{2}\,g = \left(\sum\limits_{i \,=\, 1}^{\,\infty}c_{\,i}f_{\,i}\right) \otimes \left(\sum\limits_{j \,=\, 1}^{\,\infty}d_{\,j}g_{\,j}\right) = \sum\limits_{\,i,\, j \,=\, 1}^{\,\infty}c_{\,i}\,d_{\,j}\left(\,f_{\,i} \otimes g_{\,j}\,\right)\]
On the other hand
\begin{align*}
&\left\|\,\left\{\,c_{\,i}\,\right\}_{i \,=\, 1}^{\,\infty}\,\right\|_{l^{\,2}}\,\left\|\,\left\{\,d_{\,j}\,\right\}_{j \,=\, 1}^{\,\infty}\,\right\|_{l^{\,2}} \,\leq\, C_{1}\, \left\|\,f,\, a_{\,2},\, \cdots,\, a_{\,n}\,\right\|_{\,1}\,C_{2}\, \left\|\,g,\, b_{\,2},\, \cdots,\, b_{\,n}\,\right\|_{\,2}\\
&\Rightarrow\, \left\|\,c \,\otimes\, d\,\right\|_{l^{2}} \,\leq\, C\,\left\|\,f \,\otimes\, g,\, a_{\,2} \,\otimes\, b_{\,2},\, \cdots,\, a_{\,n} \,\otimes\, b_{\,n}\,\right\|,\; \text{where}\; \,C \,=\, C_{1}\,C_{2} \,>\, 0.
\end{align*}
This completes the proof.  
\end{proof}

\begin{theorem}
If the sequence \,$\left\{\,f_{\,i} \,\otimes\, g_{\,j}\,\right\}_{i,\,j \,=\, 1}^{\,\infty}$\, is an atomic system associated to \,$\left(\,a_{\,2} \,\otimes\, b_{\,2},\, \cdots,\, a_{\,n} \,\otimes\, b_{\,n}\,\right)$\, for \,$K_{1} \,\otimes\, K_{2}$\,  in \,$H_{1} \,\otimes\, H_{2}$.\;Then \,$\left\{\,A\,f_{\,i}\,\right\}_{i \,=\, 1}^{\,\infty}$\, is an atomic system associated to \,$\left(\,a_{\,2},\, \cdots,\, a_{\,n}\,\right)$\, for \,$K_{1}$\,  in \,$H_{1}$\, and \,$\left\{\,B\,g_{\,j}\,\right\}_{j \,=\, 1}^{\,\infty}$\, is an atomic system associated to \,$\left(\,b_{\,2},\, \cdots,\, b_{\,n}\,\right)$\, for \,$K_{2}$\, in \,$H_{2}$, respectively, where \,$A$\, and \,$B$\, are constants with \,$A\,B \,=\, 1$.  
\end{theorem}

\begin{proof}
By definition (\ref{def1}), the sequence \,$\left\{\,f_{\,i} \,\otimes\, g_{\,j} \,\right\}_{i,\, j \,=\, 1}^{\,\infty}$\, is a Bessel sequence associated to \,$\left(\,a_{\,2} \,\otimes\, b_{\,2},\, \cdots,\, a_{\,n} \,\otimes\, b_{\,n}\,\right)$\, in \,$H_{1} \,\otimes\, H_{2}$, and therefore by Theorem (\ref{th3.01}), \,$\left\{\,f_{\,i}\,\right\}_{i \,=\, 1}^{\,\infty}$\, is a Bessel sequence associated to \,$\left(\,a_{\,2},\, \cdots,\, a_{\,n}\,\right)$\, in \,$H_{1}$\, and \,$\left\{\,g_{\,j}\,\right\}_{i \,=\, 1}^{\,\infty}$\, is a Bessel sequence associated to \,$\left(\,b_{\,2},\, \cdots,\, b_{\,n}\,\right)$\, in \,$H_{2}$, respectively.\;Also, for any \,$f \,\otimes\, g \,\in\, H_{F} \,\otimes\, K_{G}$, there exists \,$c \,\otimes\, d \,=\, \{\,c_{\,i}\,d_{\,j}\,\}_{i,\,j \,=\, 1}^{\,\infty}$\, in \,$l^{\,2}\,(\,\mathbb{N} \,\times\, \mathbb{N}\,)$\, such that
\[\left(\,K_{1} \,\otimes\, K_{2}\,\right)\,(\,f \,\otimes\, g\,) \,=\, \sum\limits_{i,\,j \,=\, 1}^{\,\infty}\,c_{\,i}\,d_{\,j}\,\left(\,f_{\,i} \,\otimes\, g_{\,j}\,\right) \,=\, \left(\,\sum\limits_{i \,=\, 1}^{\,\infty}\,c_{\,i}\,f_{\,i}\,\right) \,\otimes\, \left(\,\sum\limits_{j \,=\, 1}^{\,\infty}\,d_{\,j}\,g_{\,j}\,\right).\]By \,$(VI)$\, of Theorem (\ref{th1.1}), there exist constants \,$A,\,B$\, with \,$A\,B \,=\, 1$\, such that
\[K_{1}\,f \,=\, \sum\limits_{i \,=\, 1}^{\,\infty}\,c_{\,i}\,(\,A\,f_{\,i}\,)\; \;\text{and}\; \;K_{2}\,g \,=\, \sum\limits_{j \,=\, 1}^{\,\infty}\,d_{\,j}\,(\,B\,g_{\,j}\,).\]
On the other hand, for some \,$C \,>\, 0$,
\[\left\|\,c \,\otimes\, d\,\right\|_{l^{2}} \,\leq\, C\;\left\|\,f \,\otimes\, g,\, a_{\,2} \,\otimes\, b_{\,2},\, \cdots,\, a_{\,n} \,\otimes\, b_{\,n}\,\right\|, \;\text{gives}\hspace{4cm}\]
\[\left\|\,\left\{\,c_{\,i}\,\right\}_{i \,=\, 1}^{\,\infty}\,\right\|_{l^{\,2}}\left\|\,\left\{\,d_{\,j}\,\right\}_{j \,=\, 1}^{\,\infty}\,\right\|_{l^{\,2}} \,\leq\, C\left\|\,f,\, a_{\,2},\, \cdots,\, a_{\,n}\,\right\|_{\,1}\,\left\|\,g,\, b_{\,2},\, \cdots,\, b_{\,n}\,\right\|_{\,2}\;[\;\text{by (\ref{eqn1.1})}\;]\]
\[\Rightarrow \left\|\,\left\{\,c_{\,i}\,\right\}_{i \,=\, 1}^{\,\infty}\,\right\|_{l^{\,2}} \leq \dfrac{C\,\left\|\,g,\, b_{\,2},\, \cdots,\, b_{\,n}\,\right\|_{\,2}}{\left\|\,\left\{\,d_{\,j}\,\right\}_{j \,=\, 1}^{\,\infty}\,\right\|_{l^{\,2}}}\,\left\|\,f,\, a_{\,2},\, \cdots,\, a_{\,n}\,\right\|_{\,1} = C_{1}\,\left\|\,f,\, a_{\,2},\, \cdots,\, a_{\,n}\,\right\|_{\,1},\]
where \,$C_{1} \,=\, \dfrac{C\,\left\|\,g \,,\, b_{\,2} \,,\, \cdots \,,\, b_{\,n}\,\right\|_{\,2}}{\left\|\,\left\{\,d_{\,j}\,\right\}_{j \,=\, 1}^{\,\infty}\,\right\|_{l^{\,2}}} \,>\, 0$.\;Similarly, it can be shown that 
\[\left\|\,\left\{\,d_{\,j}\,\right\}_{j \,=\, 1}^{\,\infty}\,\right\|_{l^{\,2}} \,\leq\, C_{2}\,\left\|\,g,\, b_{\,2},\, \cdots,\, b_{\,n}\,\right\|_{\,2},\; \;\text{where}\; \,C_{2} \,=\, \dfrac{C\,\left\|\,f,\, a_{\,2},\, \cdots,\, a_{\,n}\,\right\|_{\,1}}{\left\|\,\left\{\,c_{\,i}\,\right\}_{i \,=\, 1}^{\,\infty}\,\right\|_{l^{\,2}}}.\]\;This completes the proof.    
\end{proof}

\begin{theorem}
Let \,$\left\{\,f_{\,i}\,\right\}_{i \,=\, 1}^{\,\infty}$\, be an atomic system associated to \,$\left(\,a_{\,2},\, \cdots,\, a_{\,n}\,\right)$\, for \,$K_{1}$\, in \,$H_{1}$\, and \,$\left\{\,g_{\,j}\,\right\}_{j \,=\, 1}^{\,\infty}$\, be an atomic system associated to \,$\left(\,b_{\,2},\, \cdots,\, b_{\,n}\,\right)$\, for \,$K_{2}$\, in \,$H_{2}$, respectively.\;Then \,$\left\{\,f_{\,i} \,\otimes\, g_{\,j}\,\right\}_{i,\,j \,=\, 1}^{\,\infty}$\, is a \,$K_{1} \,\otimes\, K_{2}$-frame associated to \,$\left(\,a_{\,2},\, \cdots,\, a_{\,n}\,\right)$.  
\end{theorem}

\begin{proof}
By Theorem (\ref{th3.01}), the sequence \,$\left\{\,f_{\,i} \,\otimes\, g_{\,j}\,\right\}_{i,\,j \,=\, 1}^{\,\infty}$\, is a Bessel sequence associated to \,$\left(\,a_{\,2} \,\otimes\, b_{\,2},\, \cdots,\, a_{\,n} \,\otimes\, b_{\,n}\,\right)$\, in \,$H_{1} \,\otimes\, H_{2}$.\;Then, for all \,$f \,\otimes\, g \,\in\, H_{F} \,\otimes\, K_{G}$, there exists \,$B \,>\, 0$\, such that 
\[\sum\limits_{i,\,j \,=\, 1}^{\,\infty}\left|\,\left<\,f \otimes g,\, f_{\,i} \otimes g_{\,j} \,|\, a_{\,2} \otimes b_{\,2},\, \cdots,\, a_{\,n} \otimes b_{\,n}\,\right>\,\right|^{\,2} \,\leq\, B \left\|\,f \otimes g,\, a_{\,2} \otimes b_{\,2},\, \cdots,\, a_{\,n} \otimes b_{\,n}\,\right\|^{\,2}\] 
Also, for any \,$f_{\,1} \,\in\, H_{F}$\, and \,$g_{\,1} \,\in\, K_{G}$, we have
\[K_{1}\,f_{\,1} \,=\, \sum\limits_{i \,=\, 1}^{\,\infty}\,c_{\,i}\,f_{\,i}\, \,\text{with}\, \,\left\|\,\left\{\,c_{\,i}\,\right\}_{i \,=\, 1}^{\,\infty}\,\right\|_{l^{\,2}} \,\leq\, C_{1}\, \left\|\,f_{\,1},\, a_{\,2},\, \cdots,\, a_{\,n}\,\right\|_{\,1},\] for some \,$C_{1} \,>\, 0$,  and
\[K_{2}\,g_{\,1} \,=\, \sum\limits_{j \,=\, 1}^{\,\infty}\,d_{\,j}\,g_{\,j}\, \,\text{with}\, \,\left\|\,\left\{\,d_{\,j}\,\right\}_{j \,=\, 1}^{\,\infty}\,\right\|_{l^{\,2}} \,\leq\, C_{2}\, \left\|\,g_{\,1} \,,\, b_{\,2} \,,\, \cdots \,,\, b_{\,n}\,\right\|_{\,2},\] for some \,$C_{2} \,>\, 0$.\;Now, for each \,$f \,\otimes\, g \,\in\, H_{F} \,\otimes\, K_{G}$, we have
\[\left\|\,\left(\,K_{1} \otimes K_{2}\,\right)^{\,\ast}(\,f \otimes g\,),\, a_{\,2} \otimes b_{\,2},\, \cdots,\, a_{\,n} \otimes b_{\,n}\,\right\|^{\,2} \,=\, \left\|\,K^{\,\ast}_{1}f \otimes K^{\,\ast}_{2}g,\, a_{\,2} \otimes b_{\,2},\, \cdots,\, a_{\,n} \otimes b_{\,n}\,\right\|^{\,2}\]
\[ \,=\, \left\|\,K^{\,\ast}_{1}\,f \,,\,  a_{\,2} \,,\, \cdots \,,\, a_{\,n}\,\right\|_{1}^{\,2}\,\left\|\,K^{\,\ast}_{2}\,g \,,\,  b_{\,2} \,,\, \cdots \,,\, b_{\,n}\,\right\|_{2}^{\,2}\; \;[\;\text{by (\ref{eqn1.1})}\;]\hspace{4cm}\]
\[ \,=\, \sup\limits_{\|\,f_{\,1},\,  a_{\,2},\, \cdots,\, a_{\,n}\,\|_{1} \,=\, 1}\left|\,\left<\,K^{\,\ast}_{1}f,\, f_{\,1} \,|\, a_{\,2},\, \cdots,\, a_{\,n}\,\right>_{1}\,\right|^{\,2}\sup\limits_{\|\,g_{\,1},\, b_{\,2},\, \cdots,\, b_{\,n}\,\|_{2} \,=\, 1}\left|\,\left<\,K^{\,\ast}_{2}g,\, g_{\,1} \,|\, b_{\,2},\, \cdots,\, b_{\,n}\,\right>_{2}\,\right|^{\,2}\]
\[ \,=\, \sup\limits_{\|\,f_{\,1},\, a_{\,2},\, \cdots,\, a_{\,n}\,\|_{1} \,=\, 1}\left|\,\left<\,f,\, K_{1}\,f_{\,1} \,|\, a_{\,2} \,,\, \cdots \,,\, a_{\,n}\,\right>_{1}\,\right|^{\,2}\sup\limits_{\|\,g_{\,1},\, b_{\,2},\, \cdots,\, b_{\,n}\,\|_{2} \,=\, 1}\,\left|\,\left<\,g,\, K_{2}\,g_{\,1} \,|\, b_{\,2},\, \cdots,\, b_{\,n}\,\right>_{2}\,\right|^{\,2}\]
\[ \,= \sup\limits_{\|\,f_{1},\, a_{2},\, \cdots,\, a_{n}\,\|_{1} = 1}\left|\left<f,\, \sum\limits_{i \,=\, 1}^{\,\infty}c_{\,i}f_{i}\,|\,a_{2},\, \cdots,\, a_{n}\right>_{1}\right|^{\,2}\sup\limits_{\|\,g_{1},\, b_{2},\, \cdots,\, b_{n}\,\|_{2} = 1}\left|\left<g,\, \sum\limits_{j \,=\, 1}^{\,\infty}d_{\,j}g_{\,j}\,|\,b_{2},\, \cdots,\, b_{n}\right>_{2}\right|^{\,2}\]
\[ = \sup\limits_{\left\|f_{1},\, a_{2},\, \cdots,\, a_{n}\right\|_{1} = 1}\left|\sum\limits_{i \,=\, 1}^{\,\infty}\overline{\,c_{\,i}}\left<f,\, f_{i}\,|\,a_{2},\, \cdots,\, a_{n}\right>_{1}\right|^{2}\sup\limits_{\|\,g_{1},\, b_{2},\, \cdots,\, b_{n}\,\|_{2} = 1}\left|\sum\limits_{j \,=\, 1}^{\,\infty}\overline{\,d_{\,j}}\left<g,\, g_{j}\,|\,b_{2},\, \cdots,\, b_{n}\right>_{2}\right|^{\,2}\]
\[\leq\, \sup\limits_{\left\|\,f_{\,1},\, a_{\,2},\, \cdots,\, a_{\,n}\,\right\|_{1} \,=\, 1}\,\left\{\,\sum\limits_{i \,=\, 1}^{\,\infty}\,\left|\,c_{\,i}\,\right|^{\,2}\,\sum\limits_{\,i \,=\, 1}^{\,\infty}\,\left|\,\left<\,f,\, f_{\,i} \,|\, a_{\,2} \,,\, \cdots \,,\, a_{\,n}\,\right>_{1}\,\right|^{\,2}\,\right\}\,\times\hspace{4cm}\]
\[\hspace{1cm}\sup\limits_{\|\,g_{\,1},\, b_{\,2},\, \cdots,\, b_{\,n}\,\|_{2} \,=\, 1}\,\left\{\,\sum\limits_{j \,=\, 1}^{\,\infty}\,\left|\,d_{\,j}\,\right|^{\,2}\,\sum\limits_{\,j \,=\, 1}^{\,\infty}\,\left|\,\left<\,g \,,\, g_{\,j} \,|\, b_{\,2},\, \cdots,\, b_{\,n}\,\right>_{2}\,\right|^{\,2}\,\right\}\]
\[\leq\, \sup\limits_{\left\|\,f_{\,1},\, a_{\,2},\, \cdots,\, a_{\,n}\,\right\|_{1} \,=\, 1}\,\left\{\,C^{\,2}_{1}\, \left\|\,f_{\,1} \,,\, a_{\,2} \,,\, \cdots \,,\, a_{\,n}\,\right\|^{\,2}_{\,1}\,\sum\limits_{i \,=\, 1}^{\,\infty}\,\left|\,\left<\,f,\, f_{\,i} \,|\, a_{\,2} \,,\, \cdots \,,\, a_{\,n}\,\right>_{1}\,\right|^{\,2}\,\right\}\,\times\hspace{1cm}\]
\[\hspace{1cm}\sup\limits_{\|\,g_{\,1},\, b_{\,2},\, \cdots,\, b_{\,n}\,\|_{2} \,=\, 1}\,\left\{\,C^{\,2}_{2}\,\left\|\,g_{\,1} \,,\, b_{\,2} \,,\, \cdots \,,\, b_{\,n}\,\right\|^{\,2}_{1}\,\sum\limits_{j \,=\, 1}^{\,\infty}\,\left|\,\left<\,g \,,\, g_{\,j} \,|\, b_{\,2},\, \cdots,\, b_{\,n}\,\right>_{2}\,\right|^{\,2}\,\right\}\]
\[=\, C^{\,2}_{1}\,C^{\,2}_{2}\;\sum\limits_{i,\,j \,=\, 1}^{\,\infty}\,\left|\,\left<\,f,\, f_{\,i} \,|\, a_{\,2} \,,\, \cdots \,,\, a_{\,n}\,\right>_{1}\,\right|^{\,2}\,\left|\,\left<\,g \,,\, g_{\,j} \,|\, b_{\,2},\, \cdots,\, b_{\,n}\,\right>_{2}\,\right|^{\,2}\hspace{3.1cm}\]
\[ \,=\, C^{\,2}_{1}\,C^{\,2}_{2}\,\sum\limits_{i,\,j \,=\, 1}^{\,\infty}\,\left|\,\left<\,f \,\otimes\, g,\, f_{\,i} \,\otimes\, g_{\,j} \,|\, a_{\,2} \,\otimes\, b_{\,2},\, \cdots,\, a_{\,n} \,\otimes\, b_{\,n}\,\right>\,\right|^{\,2}\; \;[\;\text{by (\ref{eqn1})}\;].\hspace{2cm}\]
\[\Rightarrow\dfrac{1}{C^{\,2}_{1}\,C^{\,2}_{2}}\,\left\|\,\left(\,K_{1} \,\otimes\, K_{2}\,\right)^{\,\ast}\,(\,f \,\otimes\, g\,),\, a_{\,2} \,\otimes\, b_{\,2},\, \cdots,\, a_{\,n} \,\otimes\, b_{\,n}\,\right\|^{\,2}\hspace{4cm}\]
\[\leq\, \sum\limits_{i,\,j \,=\, 1}^{\,\infty}\,\left|\,\left<\,f \,\otimes\, g,\, f_{\,i} \,\otimes\, g_{\,j} \,|\, a_{\,2} \,\otimes\, b_{\,2},\, \cdots,\, a_{\,n} \,\otimes\, b_{\,n}\,\right>\,\right|^{\,2}.\]
This completes the proof.
\end{proof}


\begin{thebibliography}{0}


\bibitem{Sadeghi}A. A. Akbar, G. Sadeghi,
\emph{Frames in 2-inner Product Spaces}, Iranian Journal of Mathematical Sciences and Informatics, Vol. 8, No. 2 (2013), pp 123-130. 

\bibitem{Janfada}Dastourian Bahram, Janfada Mohammad,
\emph{Atomic system in \,$2$-inner product spaces}, Iranian Journal of Mathematics Sciences and Informatics, vol.13, No. 1 (2018), pp 103-110.

\bibitem{Diminnie}
C. Diminnie, S. Gahler, A. White, \emph{2-inner product spaces}, Demonstratio Math. 6 (1973) 525-536.

\bibitem{L}
L. Gavruta,
\emph{Frames for operator}, Appl. Comput. Harmon. Anal. 32 (1), 139-144 (2012).

\bibitem{Mashadi}H. Gunawan, Mashadi,
\emph{On n-normed spaces}, Int. J. Math. Math. Sci., 27 (2001), 631-639. 

\bibitem{P}P. Ghosh, T. K. Samanta,
\emph{Stability of dual g-fusion frame in Hilbert spaces}, Methods of Functional Analysis and Topology, Vol. 26, no. 3, pp. 227-240.

\bibitem{Ghosh}P. Ghosh, T. K. Samanta,
\emph{Generalized atomic subspaces for operators in Hilbert spaces}, arXiv:2102.01965.

\bibitem{T} P.\,Ghosh, T.\,K Samanta,
\emph{Representation of Uniform Boundedness Principle and Hahn-Banach Theorem in linear\;$n$-normed space}, Submitted, arXiv: 2101.04555.

\bibitem{Prasenjit}P. Ghosh, T. K. Samanta,
\emph{Construction of frame relative to $n$-Hilbert space}, Submitted, arXiv: 2101.01657.

\bibitem{K}P. Ghosh, T. K. Samanta,
\emph{Generalized fusion frame in tensor product of Hilbert spaces}, arXiv:2101.06857.

\bibitem{G}P. Ghosh, T. K. Samanta,
\emph{Frame in tensor product of \,$n$-Hilbert spaces}, Submitted, arXiv: 2101.01938.

\bibitem{Misiak}A. Misiak, \emph{n-inner product spaces}, Math. Nachr., 140(1989), 299-319.

\bibitem{Folland}G.\,B. Folland,
\emph{A Course in abstract harmonic analysis}, CRC Press BOCA Raton, Florida.

\bibitem{Kadison}R.\,V Kadison and J.\,R Ringrose,
\emph{Fundamentals of the theory of operator algebras}, Vol. I, Academic Press, New York 1983.
 
\bibitem{S}S.\,Rabinson,
\emph{Hilbert space and tensor products}, Lecture notes September 8, 1997.

\bibitem{Sadri}\,V. Sadri, R. Ahmadi, A. Rahimi,
\emph{Constructions of $K$-$g$ fusion frames and their duals in Hilbert spaces}, (2018), arXiv: 1806. 03595.

\bibitem{Ahmadi}\,V. Sadri, Gh. Rahimlou, R. Ahmadi and R. Zarghami Farfar,
\emph{Generalized Fusion Frames in Hilbert Spaces}, Submitted (2018). 

\bibitem{Upender}G.\,Upender Reddy, N.\,.Gopal Reddy, \& B.\,Krishna Reddy,
\emph{Frame operator and Hilbert-Schmidt Operator in Tensor Product of Hilbert Spaces}, Journal of Dynamical Systems and Geometric Theories, 7:1, (2009), 61-70.

\bibitem{Li}Ya-Hui Wang and Yun-Zhang Li,
\emph{Tensor product dual frames},Journal of Inequalities and Applications (2019) 2019:76.

\end{thebibliography}
\end{document}